\documentclass[11pt,a4paper]{article}
\usepackage{amsfonts,amssymb,amsmath,amsthm,mathrsfs}
\newcommand{\Z}{\mathbb{Z}}
\newcommand{\R}{\mathbb{R}}
\newcommand{\C}{\mathbb{C}}

\begin{document}
\thispagestyle{empty}
\begin{center}
\textbf{{\large Differential equations and the algebra of confluent spherical functions on semisimple Lie groups.}}\\
\ \\
\ \\
\ \\
Olufemi Opeyemi {\bf OYADARE}\\
\end{center}
\ \\
\ \\
\ \\
\ \\
\ \\
\ \\
\ \\
\ \\
\ \\

\textbf{Abstract}

\begin{quote}\indent\textit{ We consider the notion of a confluent spherical function on a connected semisimple Lie group, $G,$ with finite center and of real rank $1,$ and discuss the properties and relationship of its algebra with the well-known Schwartz algebra of spherical functions on $G.$}
\end{quote}
\ \\
$\overline{\textmd{Department\;of\;Mathematics},}$
Obafemi Awolowo University,
Ile-Ife, $220005,$ Nigeria.\\
\text{E-mail: \textit{femi\_oya@yahoo.com}}
\ \\
\newpage
\indent\textbf{{\Large \S 1. Introduction}}\\

Let $G$ be a connected semisimple Lie group with finite center, a maximal compact subgroup $K$ and a Lie algebra $\mathfrak{g}$ having a \textit{Cartan decomposition,}  $$\mathfrak{g}= \mathfrak{t} \oplus \mathfrak{p}.$$ If we choose a maximal abelian subspace, $\mathfrak{a},$ of $\mathfrak{p}$ and define $A^{+} = \{ \exp tH : H \in \mathfrak{a}, t>0 \}$ then $G$ has a \textit{polar decomposition} given as $G= K \cdot cl(A^{+}) \cdot K,$ where $cl(A^{+})$ is the closure of $A^{+}.$ A function $\varphi:G\mapsto\C$ is said to be \textit{$K-$biinvariant} whenever $$\varphi(k_{1}xk_{2})=\varphi(x),\;\forall\;k_{1},k_{2}\;\in K,\;x \in G.$$ The polar decomposition of $G$ above implies that every $K-$biinvariant function on $G$ is completely determined by its restriction to $A^{+}.$ A \textit{spherical function} on $G$ is therefore a $K-$biinvariant function, $\varphi:G\mapsto\C,$ in which $\varphi(e)=1$ and which is an eigenfunction for every left-invariant differential operator on $G.$\\

An example of such a function is the Harish-Chandra \textit{(zonal) spherical function,} $\varphi_{\lambda}, \lambda \in \mathfrak{a}^{\ast}_\mathbb{\C},$ on $G.$ If we denote the restriction of $\varphi_{\lambda}$ to $A^{+}$ as $\tilde \varphi_{\lambda},$ then the following system of differential equations hold: $$\tilde q \tilde \varphi_{\lambda}= \gamma (q)(\lambda) \tilde \varphi_{\lambda},$$ where $q \in \mathfrak{Q}(\mathfrak{g}_\mathbb{\C})(:= U(\mathfrak{g}_\mathbb{\C})^{K}$ = centralizer of K in $U(\mathfrak{g}_\mathbb{\C})$), $\gamma := \gamma_{\mathfrak{g}/\mathfrak{a}}$ is the \textit{Harish-Chandra homomorphism} of $\mathfrak{Q}(\mathfrak{g}_\mathbb{\C})$ onto $U(\mathfrak{g}_\mathbb{\C})^{\mathfrak{w}},$ the $\mathfrak{w}-$ invariant subspace of $U(\mathfrak{g}_\mathbb{\C}),$ with $\mathfrak{w}$ denoting the \textit{Weyl group} of the pair $(\mathfrak{g}, \mathfrak{a}),$ $\mathfrak{t} U(\mathfrak{g}_\mathbb{\C}) \bigcap \mathfrak{Q}(\mathfrak{g}_\mathbb{\C})$ is the kernel of $\gamma$ and $\tilde q$ is the restriction of $q$ to $A^{+}.$ Since $$ \widetilde{q \cdot f}
= \tilde q \cdot \tilde f,$$ for every $f \in C^{\infty}(G//K)$ we conclude that $\tilde q$ is the \textit{radial component} of $q.$ We define $q \in \mathfrak{Q}(\mathfrak{g}_\mathbb{\C})$ to be \textit{spherical} whenever $q = \tilde q.$

The above system of differential equations have been extensively used by Harish-Chandra in the investigation of the nature of the spherical functions, $\varphi_{\lambda},$ their asymptotic expansions and their contributions to the Schwartz algebras on $G.$ The history of this investigation dated back to the $1950's$ with the two-volume work of Harish-Chandra, [$3(a.)$] and [$3(b.)$], which still attracts the strength of twenty-first century mathematicians (see [$10.$] and [$1.$]). Other functions on $G$ satisfying different interesting transformations under members of $\mathfrak{Q}(\mathfrak{g}_\mathbb{\C})$ have also been studied in the light of the approach taken by Harish-Chandra. We refer to [$5.$] and the references cited in it for further discussion.

Now if $G$ is a semisimple Lie group with real rank $1$ then it is known (see $[3a.]$) that the above system of differential equations can be replaced with $$\delta^{'}(\omega) \cdot \varphi_{\lambda}=\gamma (\omega)(\lambda) \cdot \varphi_{\lambda},$$ where $\omega$ is the \textit{Casimir operator} of $G$ and $\delta^{'}(\omega)$ denotes the radial component of the differential operator, $\delta^{'}(\omega),$ associated with $\omega.$ If we load the structure of $G,$ as a real rank $1$ semisimple Lie group, into the last equation it becomes $$(\frac{d^{2}}{dt^{2}}+ \{(p+q) \coth t+ q \tanh t \} \frac{d}{dt}) f_{\lambda}=(\lambda^{2}- \frac{(p+2q)^{2}}{4})f_{\lambda},$$ where $p=n(\alpha),$ $q=n(2\alpha),$ $f_{\lambda}(t):= \varphi_{\lambda}(\exp tH_{0})$ and $H_{0}$ is chosen in $\mathfrak{a}$ such that $\alpha(H_{0})=1$ (see $[13.],$ p. $190$ for the case of $G=SL(2, \R)$). Setting $z= -(\sinh t)^{2}$ transforms the above ordinary differential equation to the \textit{hypergeometric equation}$$(z(z-1) \frac{d^{2}}{dz^{2}} + ((a+b+1)z-c) \frac{d}{dt} + ab) g_{\lambda} = 0,$$ where $g_{\lambda}(z)=f_{\lambda}(t),$ $z<0,$ $a= \frac{p+2q+2\lambda}{4},$ $b= \frac{p+2q_2\lambda}{4}$ and $c= \frac{p+q+1}{2},$ whose solution is from here given by the \textit{Gauss hypergeometric function,} $F(a,b,c:z),$ defined as $$F(a,b,c:z)= \sum^{\infty}_{k=0} \frac{(a)_{k} (b)_{k}}{(c)_{k}} \frac{z^{k}}{k!},$$ $\mid z \mid < 1$ ([$16.$], p. $283$). It then follows that $$\varphi_{\lambda}(\exp tH_{0})= F(a,b,c:z)$$ with $z= -(\sinh t)^{2}$ and we conclude that the spherical functions on real rank $1$ semisimple Lie groups are essentially the hypergeometric function. In other words, the hypergeometric functions form the spherical functions on any real rank $1$ semisimple Lie group.

The confluent hypergeometric function is defined as $$_{1}F(a,c:z)= lim_{b\rightarrow \infty} F(a,b,c:z/b)= \sum^{\infty}_{k=0} \frac{(a)_{k}} {(c)_{k}} \frac{z^{k}}{k!}.$$ Thus replacing $z$ with $z/b$ in the hypergeometric function, $F(a,b,c:z),$ and computing the limit as $b\rightarrow \infty$ leads to its confluent, $_{1}F(a,c:z).$ Taking the same steps for the above hypergeometric equation shows that $_{1}F(a,c:z)$ satisfies the \textit{confluent differential equation} $$ z \frac{d^{2}w}{dz^{2}}+(c-z) \frac{dw}{dz}-aw=0.$$

Now since on any real rank $1$ semisimple Lie group $G,$ every spherical function is expressible, as seen above, in terms of a hypergeometric function, we refer to $_{1}F(a,c:z)$ as a \textit{confluent spherical function} on $G$ and we denote it by $\varphi^{\sigma}_{\lambda}(\exp tH_{0}).$ i.e., $$\varphi^{\sigma}_{\lambda}(\exp tH_{0})=_{1}F(a,c:z),$$ with $z=-(\sinh t)^{2}.$ It is however noted that if we replace $z$ with $z/b$ in $z=-(\sinh t)^{2},$ as in the derivation of the confluent hypergeometric equation, then we have $$t= \sinh^{-1}(i\sqrt{z/b}).$$ So that as $b\rightarrow \infty$ it implies that values of $t$ becomes very small. It then means that the relationship $$\varphi^{\sigma}_{\lambda}(\exp tH_{0})=_{1}F(a,c:z),$$ is valid only for sufficiently small values of $t.$ We conclude therefore that the spherical function, $\varphi_{\lambda}(\exp tH_{0}),$ becomes a confluent spherical function on $G$ for small values of $t.$

Our aim in this paper is, therefore, to study the function $\varphi_{\lambda}(\exp tH_{0})$ for small values of $t$ since this corresponds with the study of the confluent spherical function, $\varphi^{\sigma}_{\lambda}(\exp tH_{0})$ as explained above. In this respect we find the \textit{Stanton-Tomas expansion} of $\varphi_{\lambda}(\exp tH_{0})$ very appropriate to define the general notion of a confluent spherical function.

The paper is arranged as follows: \S 2. contains a discussion of the radial component of spherical differential operators on any $G$ of arbitrary rank, as discovered by Harish-Chandra ([$3a.$] and [$3b.$]), while the motivation for the notion of a confluent spherical function on a real rank $1$ semisimple Lie group is developed in \S3. This motivation informs our choice of the Stanton-Tomas expansion in the definition of a confluent spherical function. The algebra of these functions are then studied and related with the Schwartz algebra of spherical functions.

An insight into the study of specific confluent spherical functions on the real rank $2$ case of $Sp(2, \R),$ leading to the consideration of different kinds of Whittaker functions, is contained in Hirano, \textit{et al} $[5.].$ However the approach taken in this paper is more general than theirs and holds for \text{any} real rank $1$ semisimple Lie groups, and may be extended to higher ranks.\\

\noindent\textbf{{\Large \S 2. Radial Components of Spherical Differential Operators}}

Consider a connected real semisimple Lie group $G$ with finite center and with the Lie algebra $\mathfrak{g},$ whose complexification is denoted as $\mathfrak{g}_\mathbb{\C}.$ We can identify the members of $\mathfrak{g}$ with \textit{left-invariant vector fields} of $G$ in the following manner. For every $X \in \mathfrak{g},$ we define a map $$X\longmapsto \partial (X),$$ where $\partial (X)$ is to act on members of $C^{\infty} (G)$ by the requirement $$(\partial (X)f)(x):= \frac{d}{dt} f(x \cdot \exp tX)_{{\mid}_{t=0}}.$$ This depicts $\partial (X)$ as a first order left-invariant differential operator on $G$ associated to every $X \in \mathfrak{g}$ and which satisfies the relation $$\partial ([X,Y])=\partial (X)\partial (Y) - \partial (Y)\partial (X),$$ for $X,Y \in \mathfrak{g}.$ This outlook may be used to introduce left-invariant differential operators of any order on $G,$ by choosing more than one member of $\mathfrak{g}$ at a time. Indeed, if $X_{1},...,X_{r} \in \mathfrak{g}$ and we define the map $$X_{1} \cdots X_{r} \longmapsto \partial (X_{1} \cdots X_{r})$$ as $$\partial (X_{1} \cdots X_{r})= \partial (X_{1}) \cdots \partial(X_{r})$$ then $$\partial (X_{1} \cdots X_{r})(x)=( \frac{\partial^{r}}{\partial t_{1} \cdots \partial t_{r}} f(x \cdot \exp tX_{1} \cdots \exp tX_{r}))_{{\mid}_{(t_{1}, \cdots , t_{r})=(0, \cdots , 0)}},$$ which is a left-invariant differential operator on $G$ of order $\leq r.$ These operators are analytic and are precisely the endomorphisms of $C^{\infty}(G)$ generated by $\partial (X),$ $X \in \mathfrak{g}$ ([$13$], p. $101$). Thus if we define $D(G)= span_\mathbb{\C} \{ \partial (X) : X \in \mathfrak{g} \},$ then $D(G)$ is a subalgebra of the algebra $End_\mathbb{\C}(C^{\infty}(G)),$ of all endomorphisms on $C^{\infty}(G),$ with the identity operator as its identity element.

However, it is known that if $X_{1}, \cdots , X_{r} \in \mathfrak{g},$ the product $X_{1} \cdots X_{r}$ may not generally be a member of $\mathfrak{g},$ as may easily be verified with low-dimensional Lie algebras. Thus we should seek a \textit{gadget} in which every product, $X_{1} \cdots X_{r},$ of members of $\mathfrak{g}$ is always found, and then study the structure of the map $$\partial : X_{1} \cdots X_{r} \longmapsto \partial (X_{1} \cdots X_{r})$$ with this gadget on the foreground. With this aim in mind we consider the \textit{tensor algebra,} $T(\mathfrak{g}_\mathbb{\C}),$ of the complexification, $\mathfrak{g}_\mathbb{\C},$ of $\mathfrak{g},$ given as $$T(\mathfrak{g}_\mathbb{\C})= \mathbb{\C} \oplus \mathfrak{g}_\mathbb{\C} \oplus (\mathfrak{g}_\mathbb{\C} \otimes \mathfrak{g}_\mathbb{\C}) \oplus (\mathfrak{g}_\mathbb{\C} \otimes \mathfrak{g}_\mathbb{\C} \otimes \mathfrak{g}_\mathbb{\C}) \oplus \cdots = \bigoplus^{\infty}_{k=0} T^{k}(\mathfrak{g}_\mathbb{\C}),$$ where $T^{0}(\mathfrak{g}_\mathbb{\C}):= \mathbb{\C},$ and $T^{k}(\mathfrak{g}_\mathbb{\C}):= \mathfrak{g}_\mathbb{\C} \otimes \cdots \otimes \mathfrak{g}_\mathbb{\C}.$ $T(\mathfrak{g}_\mathbb{\C})$ is an associative algebra over $\mathbb{\C}$ with identity, and there is a natural map , $\iota,$ of $\mathfrak{g}_\mathbb{\C}$ into $T(\mathfrak{g}_\mathbb{\C})$ given by identifying $\mathfrak{g}_\mathbb{\C}$ with the first-order terms ([$1.$]). $T(\mathfrak{g}_\mathbb{\C})$ has the following universal property.\\
\ \\

$\mathbf{2.1}$\textbf{ Theorem}([$6.$], p. $644$). If $A$ is any other associative algebra over $\mathbb{\C}$ with identity and $\tau$ is a linear map of $\mathfrak{g}$ into $A,$ then there exists a unique associative algebra homomorphism $\bar \tau,$ with $\bar \tau (1)=1$ such that $\bar \tau \circ \iota = \tau.\;\Box$

We conclude, from the definitions of $\partial$ and $T(\mathfrak{g}_\mathbb{\C})$ above, that $$\partial : T(\mathfrak{g}_\mathbb{\C}) \longmapsto D(G).$$ However, as $[X,Y], XY, YX \in T(\mathfrak{g}_\mathbb{\C}),$ for every $X,Y \in \mathfrak{g}_\mathbb{\C}$ and $$\partial ([X,Y])- \partial (X) \partial (Y) - \partial (Y) \partial (X) = 0,$$ for every $X,Y \in \mathfrak{g}_\mathbb{\C},$ it would be necessary to factor, out of  $T(\mathfrak{g}_\mathbb{\C}),$ the set generated by all elements of the form $X \otimes Y - Y \otimes X - [X,Y],$ for $X,Y \in \mathfrak{g}_\mathbb{\C}.$ Indeed $$I = span_{\mathbb{\C}} \{ X \otimes Y - Y \otimes X - [X,Y] : X,Y \in \mathfrak{g}_\mathbb{\C} \}$$ is a two-sided ideal of $T(\mathfrak{g}_\mathbb{\C}),$ and we define $$U(\mathfrak{g}_\mathbb{\C}) := T(\mathfrak{g}_\mathbb{\C}) / I.$$ $U(\mathfrak{g}_\mathbb{\C})$ is also an associative algebra with identity, it contains $\mathfrak{g}_\mathbb{\C}$ and has the following universal property inherited from $T(\mathfrak{g}_\mathbb{\C}).$\\
\ \\

$\mathbf{2.2}$\textbf{ Theorem}([$6.$], p. $215$). Let $\iota$ be the canonical map of $\mathfrak{g}_\mathbb{\C}$ into $U(\mathfrak{g}_\mathbb{\C}),$ let $A$ be any complex associative algebra with identity, and let $\varphi$ be a linear mapping of $\mathfrak{g}_\mathbb{\C}$ into $A$ such that $$\varphi ([X,Y])= \varphi (X) \varphi (Y) - \varphi (Y) \varphi (X) , X,Y \in \mathfrak{g}_\mathbb{\C}.$$ Then there exists a unique algebra homomorphism  $\varphi_{0} : U(\mathfrak{g}_\mathbb{\C}) \longrightarrow A$ with $\varphi_{0}(1) = 1$ such that $ \varphi_{0} \circ \iota = \varphi.\;\Box$

The canonical map $\iota,$ in the above theorem is one-to-one ([$6.$], p. $217$) and the object $U(\mathfrak{g}_\mathbb{\C})$ is called the \textit{universal enveloping algebra} of $\mathfrak{g}_\mathbb{\C}.$ One of the most fundamental results in the theory of
$U(\mathfrak{g}_\mathbb{\C})$ is the following, which gives a concrete way of constructing it.\\
\ \\

$\mathbf{2.3}$\textbf{ Theorem}(Poincar\'{e}-Birkhoff-Witt Theorem)([$1.$], p. $32$). If $X_{1}, \cdots ,X_{n}$ is a basis of $\mathfrak{g}_\mathbb{\C}$ over $\mathbb{\C},$ then the monomials $ \iota(X_{1})^{j_{1}} \cdots \iota(X_{n})^{j_{n}},$ $j_{k} \geq 0, k=1, \cdots , n,$ form a basis of $U(\mathfrak{g}_\mathbb{\C})$ over $\mathbb{\C}.\;\Box$

The inclusion of $\iota$ in the above is not necessary since it is a one-to-one map. However, the members of $U(\mathfrak{g}_\mathbb{\C})$ may seem to be difficult to handle if we only the definition $U(\mathfrak{g}_\mathbb{\C}) := T(\mathfrak{g}_\mathbb{\C}) / I$ in mind. However recalling, from Theorem $2.2,$ that, with $A = D(G),$ $\partial : \mathfrak{g}_\mathbb{\C} \longrightarrow D(G)$ is a (natural) homomorphism such that $\partial ([X,Y])= \partial (X) \partial (Y) - \partial (Y) \partial (X) , X,Y \in \mathfrak{g}_\mathbb{\C},$ and which also extends to all of $U(\mathfrak{g}_\mathbb{\C}),$ this implies that the members of $U(\mathfrak{g}_\mathbb{\C})$ are more concrete than predicted by the Poincar\'{e}-Birkhoff-Witt Theorem. Indeed we have the following major result that gives a different outlook on $U(\mathfrak{g}_\mathbb{\C}).$\\
\ \\

$\mathbf{2.4}$\textbf{ Theorem}($[1.],$ p. $32$). The algebra homomorphism $\partial : U(\mathfrak{g}_\mathbb{\C}) \longrightarrow D(G)$ is an algebra isomorphism onto. $\Box$

The message of these theorems is that the members of $U(\mathfrak{g}_\mathbb{\C})$ are \textit{mixed derivatives,} so that $U(\mathfrak{g}_\mathbb{\C})$ is the algebra of all left-invariant differential operators on $G.$ This allows us to view $U(\mathfrak{g}_\mathbb{\C})$ as the \textit{house} of all left-invariant differential operators on $C^{\infty}(G).$ Furthermore $U(\mathfrak{g}_\mathbb{\C})$ may also be realized as the algebra of \textit{right-invariant} differential operators on $G$ via the anti-isomorphism $\partial_{r}$ given as $$(\partial_{r}(X_{1} \cdots X_{r})f)(x)=(\frac{\partial^{r}}{\partial t_{1} \cdots \partial t_{r}} f(\exp t_{1}X_{1} \cdots \exp t_{r}X_{r} \cdot x))_{{\mid}_{(t_{1}, \cdots , t_{r})=(0, \cdots , 0)}}.$$

This second realization of $U(\mathfrak{g}_\mathbb{\C})$ suggests that there are some of its members which are both left- and right-invariant. i.e., members $q \in U(\mathfrak{g}_\mathbb{\C})$ in which $qX=Xq,$ for all $X \in \mathfrak{g}_\mathbb{\C}.$ This set of members that are both left- and right-invariant is the \textit{center} of $U(\mathfrak{g}_\mathbb{\C})$ and is denoted by $\mathfrak{Z}(\mathfrak{g}_\mathbb{\C}).$ Though the algebra $\mathfrak{Z}(\mathfrak{g}_\mathbb{\C})$ is abelian and sufficient in the harmonic analysis on $G,$ we shall however consider the larger subalgebra $\mathfrak{Q}$ of $U(\mathfrak{g}_\mathbb{\C})$ defined as the centralizer of $K$ in $U(\mathfrak{g}_\mathbb{\C}).$ i.e., $$\mathfrak{Q}=\{ q \in U(\mathfrak{g}_\mathbb{\C}) : Ad(k)q=q, k \in K \}.$$ This is due to the fact that we are ultimately interested in the study of $K-$ biinvariant functions on $G.$ It is the \textit{radial component} of members of $\mathfrak{Q},$ viewed as a subalgebra of the algebra, $D(G),$ of left-invariant differential operators on $G$, that we set out to compute in this section. This is reminiscence of the classical method of finding the \textit{normal form} of an ordinary differential operator. However we need to have a generalization of the \textit{polar decomposition of matrices} to members of $G$ in order to start. We take a cue from the example of the case $G=SL(2, \R),$ where the generalization of \textit{polar coordinates} and \textit{normal form} are easily seen.

Let $G=SL(2, \R)= \{ x= \left(\begin{array}{cc} a & b \\ c & d \end{array}\right) \in GL(2, \R): ad-bc=1 \}$ with Lie algebra $\mathfrak{g}= sl(2, \R)= \{ X \in GL(2, \R) : tr(X) =0 \}$ and complexification $\mathfrak{g}_\mathbb{\C} = sl(2, \mathbb{\C}) = \{ X \in GL(2, \mathbb{\C}) : tr(X) =0 \}.$ The matrices $H= \left(\begin{array}{cc} 1 & 0 \\ 0 & -1 \end{array}\right),$ $X= \left(\begin{array}{cc} 0 & 1 \\ 0 & 0 \end{array}\right),$ and  $Y= \left(\begin{array}{cc} 0 & 0 \\ 1 & 0 \end{array}\right)$ are members of $\mathfrak{g}$ and are such that the set $\{ X-Y, H, X+Y \}$ form a basis for $\mathfrak{g}.$ Now since, for $t \in \R,$ $\exp tH= \left(\begin{array}{cc} e^{t} & 0 \\ 0 & e^{-t} \end{array}\right) =: a_{t}$ and, for $\theta \in [0, 2\pi],$ $ \exp \theta(X-Y) = \left(\begin{array}{cc} \cos \theta & \sin \theta \\ - \sin \theta & \cos \theta \end{array}\right) =: u_{\theta},$ form the closed subgroups $A$ and $K$ of $G,$ the basis of $\mathfrak{g}$ above implies that $G= K \cdot A \cdot K.$ More precisely we have $G= K \cdot cl(A^{+}) \cdot K,$ where $cl(A^{+})$ stands for the closure of $\{ a_{t} : t >0 \}.$ This is the polar decomposition which is known to generalize to any connected semisimple Lie group, with finite center. Another way to establish the polar decomposition, which easily extends to more general semisimple Lie groups, is by considering the map $$\varphi : K \times cl(A^{+}) \times K \longrightarrow G : (k_{\theta_{1}}, a_{t}, k_{\theta_{2}}) \longmapsto \varphi(k_{\theta_{1}}, a_{t}, k_{\theta_{2}})=k_{\theta_{1}}  a_{t} k_{\theta_{2}}.$$ This map is a diffeomorphism onto, and we can seek its differentials on the basis elements $X-Y, H, X+Y$ of $\mathfrak{g}.$ Indeed, since the tangent spaces to $K$ and $A$ at $1$ are $\mathfrak{t}= \{ \theta(X-Y) : \theta \in  \R \} =: \R (X-Y)$ and $\mathfrak{a}= \{ tH : t \in  \R \} =: \R H,$ then we have that $$\frac{\partial}{\partial x_{1}}= \frac{\partial}{\partial \theta_{1}},\;\frac{\partial}{\partial x_{2}}= \frac{\partial}{\partial t},\;\frac{\partial}{\partial x_{3}}= \frac{\partial}{\partial \theta_{2}},$$ in the formula $$\partial \varphi(\frac{\partial}{\partial x_{j}}_{{\mid}_{m}})= \sum^{3}_{k=1} \frac{\partial (y_{k} \circ \varphi)}{\partial x_{j}}_{{\mid}_{m}} \cdot \frac{\partial}{\partial y_{k}}_{{\mid}_{\varphi(m)}},$$ where $(x_{1},x_{2},x_{3})$ and $(y_{1},y_{2},y_{3})$ are coordinate systems about $m$ and $\varphi(m),$ respectively (cf. $[14.],$ p. $17$). Since, in this case, $m=1,$ we have $\varphi(m)=1.$ Thus, as the tangent space to $G$ at $\varphi(1)=1$ is $\mathfrak{g},$ with an orthonormal basis $\{ X-Y, H, X+Y \},$ we also have that $$\frac{\partial}{\partial y_{1}}=(X-Y),\;\frac{\partial}{\partial y_{2}}=H,\;\frac{\partial}{\partial y_{3}}=(X+Y).$$ Hence
\begin{align*}
\frac{\partial(y_{1} \circ \varphi)}{\partial x_{1}}_{{\mid}_{1}}\\
&=\frac{\partial(y_{1} \circ \varphi)}{\partial \theta_{1}}_{{\mid}_{(\theta_{1},t, \theta_{2})=(0,0,0)}}\\
&=\frac{\partial}{\partial \theta_{1}}(\exp \theta_{1}(X-Y) \cdot a_{t} \cdot k_{\theta_{2}})_{{\mid}_{(\theta_{1},t, \theta_{2})=(0,0,0)}}\\
&=\frac{\partial}{\partial \theta_{1}}(a_{t} k_{\theta_{2}}\exp \theta_{1}(X-Y)^{k_{(- \theta_{2})} a_{(- t)}})_{{\mid}_{(\theta_{1},t, \theta_{2})=(0,0,0)}}\tag{since $\exp(Z)y = y \exp(Z^{y^{-1}})$ for $Z \in \mathfrak{g}, y \in G.$}\\
&= k_{0} a_{0} \frac{d}{d \theta_{1}}(\exp \theta_{1}(X-Y)^{k_{(- \theta_{2})} a_{(- t)}})_{{\mid}_{\theta_{1}=0}}\\
&=(X-Y)^{k_{(-\theta_{2})}a_{(-t)}}\\
&=k_{(-\theta_{2})} \cdot a_{(-t)} \cdot (X-Y) \cdot a_{t} \cdot k_{\theta_{2}}\\
&=(\cosh 2t)(X-Y) + (\sin 2\theta_{2} \sinh 2t)H + (- \cos 2\theta_{2} \sinh 2t)(X+Y)\\
\end{align*}
In the same manner
\begin{align*}
\frac{\partial(y_{2} \circ \varphi)}{\partial x_{2}}_{{\mid}_{1}}=H^{k_{(-\theta_{2})}}= 0(X-Y)+(\cos 2 \theta_{2})H+(\sin 2 \theta_{2})(X+Y),\\
\frac{\partial(y_{3} \circ \varphi)}{\partial x_{3}}_{{\mid}_{1}}=(X-Y)= 1(X-Y)+ 0H + 0(X+Y).
\end{align*}
We then have the specifications of the differential, $d\varphi,$ of $\varphi$ on the basis elements, $\frac{\partial}{\partial \theta_{1}},$ $\frac{\partial}{\partial t}$ and $\frac{\partial}{\partial \theta_{2}},$ as
$$d\varphi :\left\{\begin{array}{l}
\frac{\partial}{\partial \theta_{1}} \longmapsto (\cosh 2t)(X-Y) + (\sin 2 \theta_{2} \sinh 2t)H + (- \cos 2 \theta_{2} \sinh 2t)(X+Y)\\
\ \\
\frac{\partial}{\partial t} \longmapsto 0(X-Y) + (\cos 2 \theta_{2})H + ( \sin 2 \theta_{2})(X+Y)\\
\ \\
\frac{\partial}{\partial \theta_{2}} \longmapsto 1(X-Y) + 0H + 0(X+Y)
\end{array}\right.$$

The \textit{Jacobian} of this transformation is then $ \sinh 2t,$ so that the corresponding Haar measure $,dG,$ of $G$ is also $dG= \frac{1}{2} \sinh 2t d \theta_{1} dt d \theta_{2}.$ Its inverse, $(d \varphi)^{-1},$ is then given as
$$(d \varphi)^{-1} :\left\{\begin{array}{l}
(X-Y) \longmapsto \frac{\partial}{\partial \theta_{2}}\\
H \longmapsto  \cos 2 \theta_{2} \frac{\partial}{\partial t} + \frac{ \sin 2\theta_{2}}{ \sinh 2t}(\frac{\partial}{\partial \theta_{1}}- \cosh 2t \frac{\partial}{\partial \theta_{2}})\\
(X+Y) \longmapsto \sin 2 \theta_{2} \frac{\partial}{\partial t} - \frac{ \cos 2\theta_{2}}{ \sinh 2t}(\frac{\partial}{\partial \theta_{1}}- \cosh 2t \frac{\partial}{\partial \theta_{2}})
\end{array}\right.$$
by simple substitution of the terms for $\frac{\partial}{\partial \theta_{2}},$ $\frac{\partial}{\partial t}$ and $\frac{\partial}{\partial \theta_{1}},$ respectively.

Now, since the above expression for $(d \varphi)^{-1}$ imply that $d \varphi$ is bijective everywhere on $K \times A^{+} \times K$ it follows that ($[13.],$ p. $190$) that any analytic differential operator $D$ on $G^{+}:= KA^{+}K$ gives rise to a \textit{unique} differential operator $D^{\varphi}$ on $K \times A^{+} \times K,$ called the \textit{polar form of $D,$} such that, for any $f \in C^{\infty}(G^{+}),$ we have $$(Df) \circ \varphi = D^{\varphi}(f \circ \varphi).$$ The \textit{composition with $\varphi$} in this equation means \textit{restriction to $G^{+}.$} If we now denote the restriction of $f$ to $A^{+}$ by $\tilde f,$ then the last equation above becomes $$ \widetilde {Df} = \tilde D \tilde f$$ This means that $\tilde D$ is the \textit{radial component} of the differential operator $D$ on $G^{+},$ whose existence is proved in $[3a.],$ p. $265,$ and is called \textit{spherical} whenever $D= \tilde D.$ Now since $D(G) \cong U(\mathfrak{g}_\mathbb{\C})$ it is sufficient to consider $D^{\varphi}$ for $D=Z \in \mathfrak{g}_\mathbb{\C}=sl(2, \mathbb{\C}).$ Indeed, for every $D=Z \in \mathfrak{g}_\mathbb{\C},$ we have that, $D^{\varphi}= (d \varphi)^{-1}(Z)$ and since a standard basis $ \{H^{'}, X^{'},Y^{'} \}$ of $\mathfrak{g}_\mathbb{\C}$ is given as $$H^{'}= -i(X-Y), X^{'}= \frac{1}{2}(H + i(X+Y)), Y^{'}= \frac{1}{2}(H-i(X+Y)),$$ we have that
$$H^{'{\varphi}}=(d\varphi)^{-1}(H^{'})=-i \frac{\partial}{\partial \theta_{2}},$$
$$X^{'{\varphi}}=(d\varphi)^{-1}(X^{'})= \frac{ie^{-2i \theta_{2}}}{\sinh 2t} \frac{\partial}{\partial \theta_{1}}+ e^{-2i \theta_{2}} \frac{\partial}{\partial t} - i \frac{\cosh 2t}{\sinh 2t}e^{-2i \theta_{2}} \frac{\partial}{\partial \theta_{2}},$$ and
$$Y^{'{\varphi}}=(d\varphi)^{-1}(Y^{'})= \frac{-ie^{2i \theta_{2}}}{\sinh 2t} \frac{\partial}{\partial \theta_{1}}+ e^{2i \theta_{2}} \frac{\partial}{\partial t} + i \frac{\cosh 2t}{\sinh 2t}e^{2i \theta_{2}} \frac{\partial}{\partial \theta_{2}}.$$

The Poincar\'{e}-Birkhoff-Witt Theorem above then implies that we consider the radial components of the monomials given as $$H^{'{j_{1}}} X^{'{j_{2}}} Y^{'{j_{3}}}, j_{1}, j_{2}, j_{3} \geq 0,$$ in order to exhaust the members of $U(\mathfrak{g}_\mathbb{\C}).$ It is sufficient, for a start, however to consider the center $\mathfrak{Z}(\mathfrak{g}_\mathbb{\C})$ or the \textit{centralizer,} $\mathfrak{Q}(\mathfrak{g}_\mathbb{\C}),$ of $K$ in $U(\mathfrak{g}_\mathbb{\C}),$ both of which are commutative subalgebras. Indeed, $\mathfrak{Z}(\mathfrak{g}_\mathbb{\C})$ is the commutative polynomial algebra in the single variable $$\omega =(H^{'})^{2} + 2H^{'}+4Y^{'}X^{'},$$ the \textit{Casimir operator} of $G.$ i.e., $\mathfrak{Z}(\mathfrak{g}_\mathbb{\C})= \mathbb{\C} [\omega]$ ($[8.],$ p.$195$), while $\mathfrak{Q}(\mathfrak{g}_\mathbb{\C})$ is the commutative polynomial algebra in the two variables $\omega$ and $X-Y.$ i.e., $\mathfrak{Q}(\mathfrak{g}_\mathbb{\C})= \mathbb{\C} [\omega, X-Y]$ ($[8.],$ p.$196$). Clearly $\mathfrak{Z}(\mathfrak{g}_\mathbb{\C}) \subset \mathfrak{Q}(\mathfrak{g}_\mathbb{\C}).$ We however use the \textit{normalized} casimir operator, $\omega^{'},$ given as $$\omega^{'} =(H^{'})^{2} + 2H^{'}+4Y^{'}X^{'}+1,$$ so that
\begin{align*}
\omega^{' {\varphi}}= (d\varphi)^{-1}(\omega^{'})= [\frac{1}{\sinh^{2}(2t)}(\frac{\partial^{2}}{\partial \theta^{2}_{1}}+ \frac{\partial^{2}}{\partial \theta^{2}_{2}})-2 \frac{\cosh 2t}{\sinh^{2}(2t)} \frac{\partial^{2}}{\partial \theta_{1} \partial \theta_{2}}]+[\frac{\partial^{2}}{\partial t^{2}}+2 \frac{\cosh 2t}{\sinh 2t}+1]
\end{align*}
is the restriction of $\omega^{'}$ to $G^{+},$ and may be referred to as the polar form of $\omega^{'}.$ This shows that members of $\mathfrak{Z}(\mathfrak{g}_\mathbb{\C})$ , indeed of $U(\mathfrak{g}_\mathbb{\C}),$ are essentially partial differential operators on $G^{+}.$ Hence the radial component, $\delta^{'}( \omega^{'}),$ of the normalized Casimir operator, $\omega^{'},$ which is the restriction of the above $\omega^{'}$ to $A^{+},$ is simply $$\delta^{'}( \omega^{'})= \frac{d^{2}}{dt^{2}}+2 \frac{\cosh 2t}{\sinh 2t} \frac{d}{dt}+1,$$ (see also $[7.],$ p. $73$) reducing the mixed derivatives from $\mathfrak{Z}(\mathfrak{g}_\mathbb{\C})$ to the ordinary derivatives $\frac{d^{r}}{dt^{r}}, 0\leq r\leq 2,$ on spherical functions on $G.$

The following well-known result of [$8.$], p. $199,$ explains that the eigenfunctions of $\mathfrak{Z}(\mathfrak{g}_\mathbb{\C}),$ or of $\omega^{'},$ are exactly the spherical functions on $G.$\\
\ \\

$\mathbf{2.5}$\textbf{ Theorem.} A $K-$ biinvariant $C^{\infty}$ function $f$ on $G,$ with $f(1)=1,$ is a spherical function iff $\omega^{'}= \lambda^{2}f$ for some $\lambda \in \mathbb{\C}.$ $\Box$

A detailed proof of Theorem $2.5$ is contained in [$9.$], p. $88.$ Now as every element $q \in U(sl(2, \mathbb{\C})),$ of degree $\leq r,$ may be written as $$q= \gamma H^{r}+ \sum_{l+m+n \leq r,m \leq r-1} \beta_{l,m,n}(X-Y)^{l}H^{m}(X+Y)^{n},$$ with   where $\gamma$ and $\beta_{l,m,n}$ are constants, we may generalize the above expression for $\delta^{'}( \omega^{'})$ to all members of $U(sl(2, \mathbb{\C}))$ in the following manner. Let $\mathcal{R}$ be the complex algebra of functions on $(0,\infty)$ that are generated by $(\sinh 2t)^{-1}$ and $\cosh 2t \cdot (\sinh 2t)^{-1}.$ We know that $(\frac{d}{dt}) \mathcal{R} \subset \mathcal{R},$ since the derivatives of the generators are all in $\mathcal{R}.$ The reduction of every $q \in U(sl(2, \mathbb{\C})),$ and not just of $\mathfrak{Z}(sl(2, \mathbb{\C}),$ to ordinary derivatives is established using the above expression for $q$ and the method of using $(d \varphi)^{-1}$ in the calculation of $D^{\varphi},$ for every $D=Z \in sl(2, \mathbb{\C}),$ as enumerated for $\omega^{'}$ above.\\
\ \\

$\mathbf{2.6}$\textbf{ Proposition}($[13.],$ p. $238$). If $q \in U(sl(2, \mathbb{\C}))$ is of degree $\leq r,$ then there exist $f_{0}, \cdots, f_{r-1} \in \mathbb{\C} \cdot 1 \otimes \mathcal{R}$ such that, for any $\varphi \in C^{\infty}(G//K),$ we have $$(q \cdot \varphi)=(\delta^{'}(q) \cdot \varphi)$$ on $A^{+},$ where the operator $\delta^{'}(q)$ is given as $$\delta^{'}(q)= \gamma \frac{d^{r}}{dt^{r}}+ \sum_{0 \leq j \leq r-1} f_{j} \cdot \frac{d^{j}}{dt^{j}},$$ for some constant $\gamma.\;\Box$

If $r=2$ and $q=\omega^{'} \in \mathfrak{Z}(sl(2, \mathbb{\C})),$ the conclusion of Proposition $2.6$ implies that $$\delta^{'}(\omega^{'})= \gamma \frac{d^{2}}{dt^{2}}+f_{1} \frac{d}{dt}+f_{0},$$ which is in conformity with the direct computations above, where we see that $\gamma$ is $1,$ $f_{1}(t)$ is $2 \frac{\cosh 2t}{\sinh 2t}$ and $f_{0}(t)$ is $1.$ The operator $\delta^{'}(q)$ in Propositon $2.6$ above may then be called the radial component of every $q \in U(sl(2, \mathbb{\C})).$ It would be a huge step to generalize this Proposition to every $q$ in the universal enveloping algebra, $U(\mathfrak{g}_\mathbb{\C}),$ of the complexification, $\mathfrak{g}_\mathbb{\C},$ of a real connected semisimple Lie algebra, $\mathfrak{g},$ of $G.$ To this end we extract the basic features in the above case of $G= SL(2, \R)$ as follows:

$(i.)$ Computation of the differential(s) of the \textit{Cartan decomposition map,} $K \times A^{+} \times K \longrightarrow G^{+}.$

$(ii.)$ Use this differential, in $(i.),$ to find the radial component for every member of $\mathfrak{Z}(\mathfrak{g}_ \mathbb{\C})$ or of $\mathfrak{Q}(\mathfrak{g}_ \mathbb{\C}).$

Though the programme to solve items $(i.)$ and $(ii.)$ above may not be as straightforward as we have seen for $\mathfrak{g}_ \mathbb{\C}=sl(2, \mathbb{\C}),$ the following result sets in motion the process of dealing with $(i.).$ To this end, let $\varphi : K \times A \times K \longrightarrow G$ and $G^{+}= K A^{+} K.$\\
\ \\

$\mathbf{2.7}$\textbf{ Proposition}($[2.],$ p. $125$). The map $\varphi : K \times A \times K \longrightarrow G,$ given as $\varphi(k_{1},h,k_{2}):= k_{1}hk_{2},$ is submersive on $K \times A^{+} \times K.$ In particular, $G^{+}$ is open in $G$ and $\varphi$ is an open map of $K \times A^{+} \times K$ onto $G^{+}.$ \\
\indent\textbf{\emph{Proof.}} We prove that the differential $(d \varphi)_{(k_{1},h,k_{2})} =:D$ maps $\mathfrak{t} \times \mathfrak{a} \times \mathfrak{t}$ onto $\mathfrak{g}.$ To this end, let $Z_{1},Z_{2} \in \mathfrak{t}, R \in \mathfrak{a},$ then
\begin{align*}
(d \varphi)_{(k_{1},h,k_{2})}(Z_{1},R,Z_{2}) &= D(Z_{1},R,Z_{2})\\
&= D(Z_{1},0,0) + D(0,R,0) + D(0,0,Z_{2})\\
&= (d \varphi)_{(k_{1},h,k_{2})}(Z_{1},0,0) + (d \varphi)_{(k_{1},h,k_{2})}(0,R,0) + (d \varphi)_{(k_{1},h,k_{2})}(0,0,Z_{2})\\
&= (\frac{d}{d \theta_{1}} \varphi(k_{1} \exp \theta_{1}Z_{1} \cdot h \exp 0 \cdot k_{2} \exp 0))_{\mid{_{\theta_{1}=0}}}\\
&+ (\frac{d}{dt} \varphi(k_{1} \exp 0 \cdot h \exp tR \cdot k_{2} \exp 0))_{\mid{_{t=0}}}\\
&+ (\frac{d}{d \theta_{2}} \varphi(k_{1} \exp 0 \cdot h \exp 0 \cdot k_{2} \exp \theta_{2}Z_{2}))_{\mid{_{\theta_{2}=0}}}\\
&= (\frac{d}{d \theta_{1}} \varphi(k_{1} \exp \theta_{1}Z_{1} \cdot h k_{2}))_{\mid{_{\theta_{1}=0}}} + (\frac{d}{dt} \varphi(k_{1} \cdot h \exp tR \cdot k_{2}))_{\mid{_{t=0}}}\\
&+ (\frac{d}{d \theta_{2}} \varphi(k_{1} h \cdot k_{2} \exp \theta_{2}Z_{2}))_{\mid{_{\theta_{2}=0}}}\\
&= (\frac{d}{d \theta_{1}} \varphi(k_{1} k_{2} h \cdot \exp \theta_{1}Z_{1}^{(h k_{2})^{-1}}))_{\mid{_{\theta_{1}=0}}}+ (\frac{d}{dt} \varphi(k_{1} h k_{2} \cdot \exp tR^{k^{-1}_{2}}))_{\mid{_{t=0}}}\\
&+ (\frac{d}{d \theta_{2}} \varphi(k_{1} h k_{2} \cdot \exp \theta_{2}Z_{2}))_{\mid{_{\theta_{2}=0}}}\\
&= Z^{(hk_{2})^{-1}} + R^{k^{-1}_{2}} + Z_{2}\\
&= Ad(k^{-1}_{2}) \cdot (Z^{h^{-1}}_{1}) + Ad(k^{-1}_{2}) \cdot R + Ad(k^{-1}_{2}) \cdot (Z^{k_{2}}_{2})\\
&= Ad(k^{-1}_{2}) \cdot (Z^{h^{-1}}_{1} + R + Z^{k_{2}}_{2})
\end{align*}

Hence $(Ad(k_{2}) \circ D)(Z_{1}, R, Z_{2})= Z^{h^{-1}}_{1} + R + Z^{k_{2}}_{2},$ showing that $D$ is one-to-one. The surjectivity of $D$ would
hold if we show that $\mathfrak{g} = \mathfrak{t}^{h^{-1}} + \mathfrak{a} + \mathfrak{t}$ for every $h \in A^{+}.$ Indeed, it is sufficient to verify that $\theta \mathfrak{n}$ is contained in $\mathfrak{t}^{h^{-1}} + \mathfrak{a} + \mathfrak{t},$ where $\theta$ is a \textit{Cartan involution} on $\mathfrak{g},$ $\mathfrak{n}= \sum_{\alpha \in \Delta^{+}} \mathfrak{g}_{\alpha},$ and $\Delta^{+}$ is a set of \textit{positive restricted roots} of $(\mathfrak{g},\mathfrak{a}).$ To this end, let $X \in \mathfrak{g}_{\alpha}, \alpha \in \Delta^{+},$ then $$(X+ \theta X)^{h^{-1}}= e^{- \alpha(\log h)}X + e^{\alpha(\log h)}(\theta X),$$ so that
$$ (X+ \theta X)^{h^{-1}}- e^{- \alpha(\log h)}(X+ \theta X)= e^{- \alpha(\log h)}X +  e^{\alpha(\log h)}(\theta X)- e^{- \alpha(\log h)}X - e^{- \alpha(\log h)}( \theta X)$$ $$=(e^{\alpha (\log h)}- e^{- \alpha (\log h)})(\theta X).$$ i.e.,
\begin{align*}
&\theta X = (e^{\alpha (\log h)}- e^{- \alpha (\log h)})^{-1} ((X+ \theta X)^{h^{-1}}- e^{- \alpha(\log h)}(X+ \theta X))\\
&= (e^{\alpha (\log h)}- e^{- \alpha (\log h)})^{-1} (X+ \theta X)^{h^{-1}}- (e^{\alpha (\log h)}- e^{- \alpha (\log h)})^{-1} e^{- \alpha(\log h)}(X+ \theta X) \cdots (*)\\
& \in \mathfrak{t}^{h^{-1}}+\mathfrak{t} \tag{since $(X+ \theta X) \in \mathfrak{t}$}\\
& \subset \mathfrak{t}^{h^{-1}} + \mathfrak{a} + \mathfrak{t}.
\end{align*}
This ends the proof as expected.$\;\Box$

We see that $(d \varphi)_{(1,h,1)}= Z^{h^{-1}}_{1} + R + Z_{2},$ in anticipation of its use on the $K-$ biinvariant functions on $G.$ The proof of Proposition $2.7$ above gives a formula for the first differential, $(d \varphi)_{(k_{1},h,k_{2})},$ of $\varphi.$ However since we are ultimately interested in the \textit{radial component} of an arbitrary $C^{\infty}$ \textit{spherical} differential operator on $G,$ which may have second, third, and higher derivatives, we compute higher order derivatives of $\varphi$ to give the full differential, which we shall denote by $(d \varphi)^{\infty}_{(k_{1},h,k_{2})}.$ This is mainly because the property of $\mathfrak{Z}(sl(2, \mathbb{\C})),$ as a polynomial algebra in the variable $\omega^{'},$ has not be found generalizeable to arbitrary semisimple $G$ and $\mathfrak{g}.$ Indeed we have the following.\\
\ \\

$\mathbf{2.8}$\textbf{ Proposition}($[2.],$ p. $127$). The full differential, $(d \varphi)^{\infty}_{(k_{1},h,k_{2})},$ of the map $\varphi$ in Proposition $2.7$ is given on $U(\mathfrak{t}_\mathbb{\C}) \oplus U(\mathfrak{a}_{_\mathbb{\C}}) \oplus U(\mathfrak{t}_{_\mathbb{\C}})$ as $$(d \varphi)^{\infty}_{(k_{1},h,k_{2})}(\xi_{1} \otimes u \otimes \xi_{2})= \xi^{(hk_{2})^{-1}}_{1}u^{k^{-1}_{2}}\xi_{2},$$ where $\xi_{1},\xi_{2} \in U(\mathfrak{t}_\mathbb{\C})$ and $u \in U(\mathfrak{a}_\mathbb{\C}).$ In particular $$(d \varphi)^{\infty}_{(1,h,1)}(\xi_{1} \otimes u \otimes \xi_{2})= \xi^{h^{-1}}_{1}u \xi_{2}.\;\Box$$

Just as in Proposition $2.7,$ the map $$D_{h}:U(\mathfrak{t}_\mathbb{\C}) \oplus U(\mathfrak{a}_{_\mathbb{\C}}) \oplus U(\mathfrak{t}_{_\mathbb{\C}}) \longrightarrow U(\mathfrak{g}_\mathbb{\C})$$ defined as $D_{h}= (d \varphi)^{\infty}_{(1,h,1)}$ is surjective. Thus for every $q \in U(\mathfrak{g}_\mathbb{\C})$ there exists $\tau_{h} \in U(\mathfrak{t}_\mathbb{\C}) \oplus U(\mathfrak{a}_{_\mathbb{\C}}) \oplus U(\mathfrak{t}_{_\mathbb{\C}})$ such that $D_{h}(\tau_{h})=q.$ If we assume that $\tau_{h}$ depends smoothly on $h,$ then the map $h \longmapsto \tau_{h}$ leads to a differential operator on $K \times A^{+} \times K$ which, at the points of $(1) \times A^{+} \times (1),$ is simply $q$ \textit{expressed in polar coordinates.} To the find the formula for $\tau_{h}$ we proceed as follows.

Let $q \longmapsto t(q)$ be the projection of $ U(\mathfrak{g}_\mathbb{\C})$ onto $U(\mathfrak{a}_\mathbb{\C}) U(\mathfrak{t}_\mathbb{\C}),$ which corresponds to the direct sum $U(\mathfrak{g}_\mathbb{\C})=  U(\mathfrak{a}_\mathbb{\C}) U(\mathfrak{t}_\mathbb{\C}) \oplus \theta (\mathfrak{n}) U(\mathfrak{g}_\mathbb{\C}),$ then $deg(t(q)) \leq deg(q).$ In view of the $f_{j}'s$ in Proposition $2.6,$ we define $f_{\alpha}$ and $g_{\alpha}, \alpha \in \Delta^{+},$ on $A^{+}$ by $f_{\alpha}= (\xi_{\alpha}-\xi_{- \alpha})^{-1}$ and $g_{\alpha}= \xi_{- \alpha}(\xi_{\alpha}-\xi_{- \alpha})^{-1},$ respectively, where $\xi_{\lambda}= e^{\lambda \circ \log}, \lambda \in \mathfrak{a}_\mathbb{\C}.$ Also, let $\mathcal{R}_{0}$ be the algebra with unit generated over, $\mathbb{\C},$ by the $f_{\alpha}$  and $g_{\alpha},$ and, for any integer $d \geq 1,$ let $\mathcal{R}_{0,d}$ be the linear span of the monomials in these generators of degree $d.$ Now if we put $$\mathcal{R}^{+}_{0}= \sum_{d \geq 1} \mathcal{R}_{0,d},$$ then for every $q \in U(\mathfrak{g}_\mathbb{\C})$ of degree $m,$ there exist $\xi_{i},\xi^{'}_{i} \in U(\mathfrak{t}_\mathbb{\C}), u_{i} \in U(\mathfrak{a}_\mathbb{\C})$ and $\varphi_{i} \in \mathcal{R}^{+}_{0}, 1 \leq i \leq n$ such that $$q= t(q)+ \sum_{1 \leq i \leq n} \varphi_{i}(h)\xi^{h^{-1}}_{i}u_{i}\xi^{'}_{i},\;\;h \in A^{+},$$ and we may take $$\tau_{h}= 1 \otimes t(q)+ \sum_{1 \leq i \leq n} \varphi_{i}(h)\xi_{i} \otimes u_{i} \otimes \xi^{'}_{i}, h \in A^{+},\;\;q \in U(\mathfrak{a}_\mathbb{\C}).$$ This is the appropriate generalization of the expression for $q$ in $U(sl(2,\mathbb{\C}))$ of degree $\leq r$ as given after Theorem $2.5$ above. This direct comparison with the case of $sl(2, \mathbb{\C})$ implies that we need to use the general expression for $q$ above to seek the generalization of Proposition $2.6$ to all members of $U(\mathfrak{g}_\mathbb{\C})$ for any semisimple Lie algebra $\mathfrak{g}.$ It is sufficient, in our present case, to seek this generalization to all members of $\mathfrak{Q}(\mathfrak{g}_\mathbb{\C})$ as we now do next. First a little preparation.

Take $\beta_{\mathfrak{n}}: U(\mathfrak{g}_\mathbb{\C}) \longrightarrow U(\mathfrak{a}_\mathbb{\C})$ be the projection corresponding to the direct sum $U(\mathfrak{g}_\mathbb{\C})= U(\mathfrak{a}_\mathbb{\C}) \oplus (\mathfrak{t} U(\mathfrak{g}_\mathbb{\C}) + U(\mathfrak{g}_\mathbb{\C}) \mathfrak{n}),$ having $\mathfrak{Q} \cap (\mathfrak{t} U(\mathfrak{g}_\mathbb{\C})$ as its kernel ($[3(a.)],$ p. $260$) and define a map $\gamma_{\mathfrak{n}}: \mathfrak{Q}(\mathfrak{g}_{\mathbb{\C}}) \longrightarrow U(\mathfrak{a}_\mathbb{\C})$ by the specification $$\gamma_{\mathfrak{n}}(q)(\lambda)= \beta_{\mathfrak{n}}(q)(\lambda - \rho), \lambda \in \mathfrak{a}^{*}_{\mathbb{\C}},q \in \mathfrak{Q}(\mathfrak{g}_{\mathbb{\C}}), \rho = \frac{1}{2} \sum_{\alpha \in \Delta^{+}} dim(\mathfrak{g}_{\alpha}) \cdot \alpha.$$ $\gamma_{\mathfrak{n}}$ is a homomorphism ($[3(a.)],$ p. $260$), is independent of the choice of $\mathfrak{n}$ and is called the \textit{Harish-Chandra homomorphism.} We  denote it simply as $\gamma$ which, for $\mathfrak{g}= sl(2, \R),$ is given on $\omega^{'}= (H^{'})^{2} + 2H^{'} + 4Y^{'}X^{'}$ as $\gamma(\omega^{'})=H^{2} \longmapsto \frac{d^{2}}{dt^{2}}.$ (see [$7.$], p. $51$ and use the isomorphism in Theorem $2.4$) We state the major result of this section.\\
\ \\

$\mathbf{2.9}$\textbf{ Theorem}($[2.],$ p. $129$ and $[3(a.)],$ p. $267$). Given any analytic spherical differential operator $E$ on $G^{+}$ there is a unique analytic differential operator $\tilde E$ on $A^{+},$ the radial component of $E,$ such that $\widetilde{Ef}= \tilde E \tilde f, f \in C^{\infty}(G//K).$ The map $E \longmapsto \tilde E$ is a homomorphism that does not increase degree. If $E=q \in \mathfrak{Q}(\mathfrak{g}_{\mathbb{\C}}),$ then the radial component, written as $\delta^{'}(q),$ is given as $$\delta^{'}(q)= e^{- \rho} \cdot \gamma(q) \cdot e^{\rho}+ \sum_{1 \leq i \leq n} \varphi_{i}u_{i},$$ where $u_{i} \in U(\mathfrak{g}_{\mathbb{\C}}), \varphi_{i} \in \mathcal{R}^{+}_{0}$ and $deg(u_{i}) < deg(q).\;\Box$

It is clear that Theorem $2.9$ generalizes the assertions of Proposition $2.6,$ at least to all members of $\mathfrak{Q}(\mathfrak{g}_{\mathbb{\C}}),$ and it sets the stage for analysis of the differential equations satisfied by spherical functions on $G.$ Indeed, by Theorem $2.5$ we have that $ \omega^{'} \varphi = \lambda^{2} \varphi$ which, when combined with Proposition $2.9,$ (see also Lemma $23$ of $[3(a.)]$) implies that $\delta^{'}(\omega^{'}) \varphi= \lambda^{2} \varphi$ where, according to Theorem $2.9,$ $$\delta^{'}(\omega^{'})= e^{- \rho} \cdot \gamma(\omega^{'}) \cdot e^{\rho}+ \sum_{1 \leq i \leq n} \varphi_{i}u_{i},$$ and $\omega^{'}$ is the \textit{normalized} Casimir operator of $U(\mathfrak{g}_{\mathbb{\C}}).$ Even though the relation $\mathfrak{Z}(sl(2, \mathbb{\C}))= \mathbb{\C}[\omega^{'}]$ does not generalize to arbitrary $\mathfrak{g}_{\mathbb{\C}}$ the sufficiency of considering the differential equations $$\delta^{'}(\omega^{'}) \varphi= \lambda^{2} \varphi \;\;(= \gamma_{\mathfrak{n}}(\omega^{'})(\lambda) \cdot \varphi)$$ in our study of spherical functions on $G$ may be justified as follows.\\
\ \\

$\mathbf{2.10}$\textbf{ Theorem}($[2.],$ p. $145$). If $\varphi_{\lambda}, \lambda \in \mathfrak{a}^{*}_{\mathbb{\C}}$ is an eigenfunction of $\delta^{'}(\omega^{'}),$ then it is also an eigenfunction of $\delta^{'}(q), q \in \mathfrak{Q}(\mathfrak{g}_{\mathbb{\C}}),$ with the same eigenvalue.\\
\indent\textbf{\emph{Proof.}} We first show that the differential operators, $\delta^{'}(q), q \in \mathfrak{Q}(\mathfrak{g}_{\mathbb{\C}}),$ commute with each other. Indeed, for any $q_{1},q_{2} \in \mathfrak{Q}(\mathfrak{g}_{\mathbb{\C}}),$ and any $f \in C^{\infty}(G^{+}//K),$ the commutativity of $\mathfrak{Q}(\mathfrak{g}_{\mathbb{\C}})$ implies that $q_{1}q_{2}f=q_{2}q_{1}f.$ Hence $$\delta^{'}(q_{1}q_{2}) \tilde f = \delta^{'}(q_{2}q_{1}) \tilde f,$$ where$ \tilde f= f_{\mid_{A^{+}}}.$ Therefore $$\delta^{'}(q_{1}) \delta^{'}(q_{2}) = \delta^{'}(q_{2}) \delta^{'}(q_{1}) \cdots (*)$$

Now let $\varphi_{\lambda}, \lambda \in \mathfrak{a}^{*}_{\mathbb{\C}},$ be an eigenfunction of $\delta^{'}(\omega^{'}).$ i.e., for some $\lambda \in \mathbb{\C},$ we have $\delta^{'}(\omega^{'})\varphi_{\lambda}= \gamma(\omega^{'})(\lambda) \cdot \varphi_{\lambda}.$ Thus $$\delta^{'}(\omega^{'})(\delta^{'}(q) \varphi_{\lambda})= \delta^{'}(q)(\delta^{'}(\omega^{'}) \varphi_{\lambda})= \gamma(\omega^{'})(\lambda)(\delta^{'}(q) \varphi_{\lambda}),$$ meaning that $\delta^{'}(q) \varphi_{\lambda}$ is an eigenfunction of $\delta^{'}(\omega^{'}).$ However, by the uniqueness of Theorem $2.9,$ we must have that $\delta^{'}(q) \varphi_{\lambda}$ is a constant multiple of $\varphi_{\lambda,}$ as required. $\Box$

The above result explains that it is sufficient to consider the operator $$\delta^{'}(\omega^{'})= \tilde \omega^{'}= \frac{d^{2}}{dt^{2}}+ 2 \frac{\cosh 2t}{\sinh 2t} \frac{d}{dt}+1$$ in the study of spherical functions, $\varphi_{\lambda},$ on $SL(2, \R).$ In this case we have $\gamma(\omega^{'})(\lambda)= \lambda^{2}$ so that the equation $$\delta^{'}(\omega^{'}) \varphi_{\lambda}= \gamma_{\mathfrak{n}}(\omega^{'})(\lambda) \cdot \varphi_{\lambda}$$ becomes $$\frac{d^{2} \varphi_{\lambda}}{dt^{2}}+ 2 \frac{\cosh 2t}{\sinh 2t} \frac{d \varphi_{\lambda}}{dt}+ (1-\lambda^{2}) \varphi_{\lambda} = 0 \cdots (\Upsilon)$$ i.e., $$\sinh^{2}(2t) \frac{d^{2} \varphi_{\lambda}}{dt^{2}}+ 2 \sinh 2t \cosh 2t \frac{d \varphi_{\lambda}}{dt}+ \sinh^{2}(2t)(1-\lambda^{2}) \varphi_{\lambda} = 0 $$ Now setting $z= \cosh 2t,$ and defining $\varphi_{\lambda}(t)$ as $\Phi_{\lambda}(z),$ we would have $\frac{d \varphi_{\lambda}}{dt}= (2 \sinh 2t) \frac{d \Phi_{\lambda}}{dz}$ and $\frac{d^{2} \varphi_{\lambda}}{dt^{2}}= 4 \sinh^{2}(2t) \frac{d^{2} \Phi_{\lambda}}{dz^{2}}+4 \cosh 2t \frac{d \Phi_{\lambda}}{dz}$ so that the last differential equation above transforms to $$\sinh^{4}(2t) \frac{d^{2} \Phi_{\lambda}}{dz^{2}}+ 2 \sinh^{2}(2t) \cosh (2t) \frac{d \Phi_{\lambda}}{dz} + \frac{1}{4} \sinh^{2}(2t)(1- \lambda^{2}) \Phi_{\lambda}= 0.$$ We finally have, with $\sinh^{2}(2t)=-(1-z^{2}),$ $$(1-z^{2}) \frac{d^{2} \Phi_{\lambda}}{dz^{2}}- 2z \frac{d \Phi_{\lambda}}{dz} + \frac{(\lambda^{2}-1)}{4} \Phi_{\lambda}= 0.$$

This is the well-known \textit{Legendre equation.}

{\bf This comfirms that the spherical functions on $G= SL(2,\R)$ are essentially the Legendre functions as enunciated in [$4.$], pp. $405-407.$}

It is a well-known fact in the general theory of ordinary differential equations that we can now consider the associated \textit{confluent Legendre functions} on $G= SL(2,\R)$ in this context, and since every spherical function on $G= SL(2,\R)$ is a Legendre function, we may refer to the confluent Legendre functions on $G= SL(2,\R)$ as \textit{confluent spherical functions.} This is the motivation for the next section where this outlook is generated to all semisimple Lie groups with real rank $1.$\\

\noindent\textbf{{\Large \S 3. Reduction to the real rank 1 case.}}

It is appropriate, from Theorem $2.10,$ to find the general Casimir operator for the semisimple Lie group $G$ in order to get the generalization
of the Legendre equation of the group $G= SL(2,\R)$ as already seen above. To this end let $J$ be the two-sided proper ideal of $U(\mathfrak{g}_{\mathbb{\C}}),$ generated over $\mathbb{\C}$ by elements of the form $X \otimes Y - Y \otimes X,$ where $X,Y \in \mathfrak{g}_{\mathbb{\C}}.$ The quotient $T(\mathfrak{g}_{\mathbb{\C}})/J$ is the \textit{symmetric algebra} of $\mathfrak{g}_{\mathbb{\C}},$ denoted as $\mathcal{S}(\mathfrak{g}_{\mathbb{\C}}).$ Clearly if $\mathfrak{g}_{\mathbb{\C}}$ is abelian, then $U(\mathfrak{g}_{\mathbb{\C}})=\mathcal{S}(\mathfrak{g}_{\mathbb{\C}})$ Even if $\mathfrak{g}_{\mathbb{\C}}$ is non abelian, so that we only have $U(\mathfrak{g}_{\mathbb{\C}}) \supset \mathcal{S}(\mathfrak{g}_{\mathbb{\C}})$ in general, there is a map, $\lambda : \mathcal{S}(\mathfrak{g}_{\mathbb{\C}}) \longrightarrow U(\mathfrak{g}_{\mathbb{\C}}),$ called the \textit{Harish-Chandra symmetrization map} given as $$\lambda(X_{1} \cdots X_{r})= \frac{1}{r!} \sum_{\sigma}X_{\sigma(1)} \cdots X_{\sigma(r)},$$ where $X_{1}, \cdots X_{r} \in \mathfrak{g}_{\mathbb{\C}}$ and $\sigma$ runs over the set of all permutations of the set $\{1, \cdots ,r \}.$ This is a linear isomorphism and it may be shown that, for every $x \in G,$ we have $\lambda \circ Ad(x)= Ad(x) \circ \lambda,$ where $Ad(x)$ is viewed as a map from $\mathcal{S}(\mathfrak{g}_{\mathbb{\C}})$ into $U(\mathfrak{g}_{\mathbb{\C}})$ (By an adaptation of Theorem $2.2$ to $\mathcal{S}(\mathfrak{g}_{\mathbb{\C}})$)([$4.$], p. $393$). We the have the following.\\
\ \\

$\mathbf{3.1}$\textbf{ Theorem}($[13.],$ p. $103$). If $I(\mathfrak{g}_{\mathbb{\C}})$ is the subset of members of $\mathcal{S}(\mathfrak{g}_{\mathbb{\C}})$ which are $Ad(G)-$ invariant, then $\lambda : I(\mathfrak{g}_{\mathbb{\C}}) \longrightarrow \mathfrak{Z}(\mathfrak{g}_{\mathbb{\C}})$ is a linear isomorphism.\\
\indent\textbf{\emph{Proof.}} This is a direct consequence of the relation $\lambda \circ Ad(x)= Ad(x) \circ \lambda,$ for $x \in G.\;\Box$

The above result allows us to use a basis of $I(\mathfrak{g}_{\mathbb{\C}})$ in the construction of a basis of $\mathfrak{Z}(\mathfrak{g}_{\mathbb{\C}}).$ Since it is known that $\mathcal{S}(\mathfrak{g}_{\mathbb{\C}})$ is essentially the polynomial algebra on $\mathfrak{g}_{\mathbb{\C}},$ we may introduce the \textit{Casimir polynomial,} $\xi,$ on $\mathfrak{g}_{\mathbb{\C}}$ given as $\xi(Z)= tr(ad Z)^{2}.$ It follows that $\xi \in \mathcal{S}(\mathfrak{g}_{\mathbb{\C}})$ and hence we have $\lambda(\xi) \in \mathfrak{Z}(\mathfrak{g}_{\mathbb{\C}}).$ If we define $\tilde \xi$ on $\mathfrak{g}_{\mathbb{\C}}$ by the requirement $$\langle \tilde \xi, Z \rangle= \xi(Z), Z \in \mathfrak{g}_{\mathbb{\C}},$$ then we shall refer $\lambda(\tilde \xi)$ as the Casimir operator and denote it by $\omega.$ The situation for $\mathfrak{g}= sl(2, \R)$ may be used to justify these terms. Indeed, the Casimir polynomial in this example is $H^{2}+4YX$ and, hence, $\lambda(H^{2}+4YX)=\lambda(HH)+4 \lambda(YX)= \frac{1}{2!}(HH+HH)+ 4[\frac{1}{2!}(YX+XY)]=H^{2}+2YX+2XY=H^{2}+2YX+(2H+2YX)=H^{2}+2H+4YX,$ which when normalized gives exactly $\omega^{'}.$ In the general case we have the following, where we denote $B$ as the Cartan-Killing form on $\mathfrak{g}_{\mathbb{\C}} \times \mathfrak{g}_{\mathbb{\C}}.$\\
\ \\

$\mathbf{3.2}$\textbf{ Theorem}($[12.],$ p. $217$). The Casimir operator $\omega$ belongs to the center, $\mathfrak{Z}(\mathfrak{g}_{\mathbb{\C}})$ of $\mathfrak{Y}(\mathfrak{g}_{\mathbb{\C}}).$ If $\{X_{1}, \cdots, X_{m} \}$ is a basis for $\mathfrak{g}_{\mathbb{\C}}$ and $\{X^{1}, \cdots, X^{m} \}$ is the dual basis defined by $B(X_{i},X^{j})=\delta_{ij},$ then $$\omega= \sum_{1 \leq i \leq m}X_{i}X^{i}.$$
\indent\textbf{\emph{Proof.}} The first statement holds from the linear isomorphism in Theorem $3.1.$ The dual basis, $\{X^{1}, \cdots, X^{m} \},$ exists since $B$ is a non-singular symmetric bilinear form. Now for every $X \in \mathfrak{g}_{\mathbb{\C}}$ we have that $X= \sum_{1 \leq j \leq m}B(X,X^{j})X_{j}$ so that $$(adX)^{2}= \sum_{1 \leq r,s \leq m}B(X,X^{r})B(X,X^{s})adX_{r}adX_{s}.$$ Therefore $$\xi(X)=tr(adX)^{2}= \sum_{1 \leq r,s \leq m} B(X,X^{r})B(X,X^{s})B(X_{r},X_{s}).$$ Hence $\tilde \xi= \sum_{1 \leq r,s \leq m}B(X_{r},X_{s})X^{r}X^{s}.$ We then have that $\omega= \lambda(\tilde \xi)\\
= \sum_{1 \leq r,s \leq m}B(X_{r}X_{s})X^{r}X^{s}= \sum_{1 \leq s \leq m}(\sum_{1 \leq r \leq m}B(X_{r},X_{s})X^{r})X^{s}= \sum_{1 \leq s \leq m}X_{s}X^{s},$ as expected. $\Box$

We need to now compute the expression for the constant coefficient differential operator, $\gamma(\omega),$ for $\omega$ as in Theorem $3.2$ above, in anticipation of its use in the formula for $\delta^{'}(\omega)$ contained in Theorem $2.9.$ This would require using a specific basis of $\mathfrak{g}$ according to the decomposition $\mathfrak{g}= \mathfrak{m} \oplus \mathfrak{a} \oplus \sum_{\lambda \in \Delta} \mathfrak{g}_{\lambda},$ as we now proceed next.

Let $\{Z_{i} \}_{1 \leq i \leq m}$ be a basis of $\mathfrak{m}$ such that $B(Z_{i},Z_{j})= - \delta_{ij},$ and let $\{H_{i} \}_{1 \leq i \leq a}$ be a basis of $\mathfrak{a}$ such that $B(H_{i},H_{j})= \delta_{ij}.$ This is possible since $B$ is negative-definite on $\mathfrak{t}$ (hence on $\mathfrak{m}$) and positive-definite on $\mathfrak{p}$ (hence on $\mathfrak{a}$). Also since, for any two roots $\alpha, \beta$ of $(\mathfrak{g},\mathfrak{a}),$ the root-space $\mathfrak{g}_{\alpha}$ is orthogonal to $\mathfrak{g}_{\beta}$ whenever $\alpha \neq \beta,$ and since $B$ is non-degenerate when restricted to $\mathfrak{g}_{\alpha} \times \mathfrak{g}_{- \alpha}$ (see [$4.$], p. $141$) we may select a basis $\{X_{\alpha,i} \}_{1 \leq i \leq dim(\mathfrak{g}_{\alpha})},$ of $\mathfrak{g}_{\alpha},$ such that $B(X_{\alpha,i},X_{- \alpha,i})= \delta_{ij},$ where $\theta X_{\alpha,i}= X_{- \alpha,i}$ for a Cartan involution, $\theta,$ on $\mathfrak{g},$ for every root $\alpha.$ It is clear, from the assertion of Theorem $3.2,$ that $\omega_{\mid_{\mathfrak{m}}}= \sum_{1 \leq i,j \leq m}Z_{i}Z^{j}=-(Z^{2}_{1}+ \cdots +Z^{2}_{m}),$ $\omega_{\mid_{\mathfrak{a}}}= \sum_{1 \leq i,j \leq a}H_{i}H^{j}=H^{2}_{1}+ \cdots +H^{2}_{a}$ and $\omega_{\mid_{\mathfrak{g}_{\alpha}}}= \sum_{1 \leq i,j \leq dim(\mathfrak{g}_{\alpha})}X_{\alpha,i}X^{\alpha,j}= \sum_{1 \leq i,j \leq dim(\mathfrak{g}_{\alpha})}(X_{\alpha,i}X_{- \alpha,i}+ X_{- \alpha,i}X_{\alpha,i}).$ Hence the direct sum $\mathfrak{g}= \mathfrak{m} \oplus \mathfrak{a} \oplus \sum_{\lambda \in \Delta} \mathfrak{g}_{\lambda}$ now implies that $$\omega=-(Z^{2}_{1}+ \cdots +Z^{2}_{m})+ (H^{2}_{1}+ \cdots +H^{2}_{a})+ \sum_{1 \leq i \leq n(\alpha), \alpha > 0}(X_{\alpha,i}X_{- \alpha,i}+ X_{- \alpha,i}X_{\alpha,i}).$$ The present form of the Casimir operator in theses bases may now be used to compute the constant coefficient differential operator, $\gamma(\omega).$\\
\ \\

$\mathbf{3.3}$\textbf{ Lemma}($[2.],$ p. $94$). In the above form of $\omega,$ we have that $$\gamma(\omega)= \sum_{1 \leq i \leq a}H^{2}_{i}-B(H_{\rho},H_{\rho})$$ where $H_{\rho}= \frac{1}{2} \sum_{\alpha > 0}n(\alpha)H_{\alpha}$ and $H_{\alpha}$ is uniquely defined by the requirement, $\alpha(H)=B(H_{\alpha},H)$ for all $H \in \mathfrak{a}.$

\indent\textbf{\emph{Proof.}} We know that $X_{\alpha,i}X_{- \alpha,i}=X_{- \alpha,i}X_{\alpha,i}+[X_{\alpha,i}, X_{- \alpha,i}].$ Now if $X \in \mathfrak{g}_{\alpha}$ and $X^{'} \in \mathfrak{g}_{- \alpha},$ then $[X,X^{'}] \in \mathfrak{g}_{\alpha = 0} \subset \mathfrak{m} \oplus \mathfrak{a},$ and, for every $H \in \mathfrak{a},$ we always have $B(H,[X,X^{'}])=\alpha(H)B(X,X^{'}).$ These sum up to give $[X,X^{'}] \equiv B(X,X^{'})H_{\alpha}(mod \mathfrak{m}),$ for every $X \in \mathfrak{g}_{\alpha}, X^{'} \in \mathfrak{g}_{- \alpha}.$ In particular $$[X_{\alpha,i}, X_{- \alpha,i}] \equiv B(X_{\alpha,i}, X_{- \alpha,i})H_{\alpha}(mod \mathfrak{m})= H_{\alpha}(mod \mathfrak{m}).$$ It then follows that
\begin{align*}
\omega &= - \sum^{m}_{k=1}Z^{2}_{k}+ \sum^{a}_{i=1} H^{2}_{i}+ \sum_{1 \leq i \leq n(\alpha), \alpha > 0}(X_{\alpha,i}X_{- \alpha,i}+ X_{- \alpha,i}X_{\alpha,i})\\
&= - \sum^{m}_{k=1}Z^{2}_{k}+ \sum^{a}_{i=1} H^{2}_{i}+ \sum_{1 \leq i \leq n(\alpha), \alpha > 0}(2 X_{- \alpha,i}X_{\alpha,i}+ [X_{\alpha,i}, X_{- \alpha,i}])\\
& \equiv - \sum^{m}_{k=1}Z^{2}_{k}+ \sum^{a}_{i=1} H^{2}_{i}+ \sum_{1 \leq i \leq n(\alpha), \alpha > 0}(2 X_{- \alpha,i}X_{\alpha,i}+ H_{\alpha}(mod \mathfrak{m})).\\
\end{align*}
Now using the fact that $\mathfrak{Q}(\mathfrak{g}_{\mathbb{\C}}) \cap (\mathfrak{t} U(\mathfrak{g}_{\mathbb{\C}})),$ which contains $\mathfrak{m},$ is the kernel of $\beta_{\mathfrak{n}},$ we have $$\beta_{\mathfrak{n}}= -0 + \sum^{a}_{i=1} H^{2}_{i} + \sum_{\alpha > 0}n(\alpha)H_{\alpha}.$$ Hence
\begin{align*}
\gamma(\omega)(\lambda) &= \beta_{\mathfrak{n}}(\omega)(\lambda - \rho)\\
&= \beta_{\mathfrak{n}}(\omega)(\lambda)- \beta_{\mathfrak{n}}(\omega)(\rho)\\
&= (\sum^{a}_{i=1} H^{2}_{i} + \sum_{\alpha > 0}n(\alpha)H_{\alpha})- (2 \sum_{1 \leq i \leq a} \rho(H_{i})H_{i}- \sum_{1 \leq i \leq a} \rho(H_{i})^{2} + \sum_{\alpha > 0}n(\alpha) \rho(H_{\alpha})).
\end{align*}
i.e., $$\gamma(\omega)= \sum^{a}_{i=1} H^{2}_{i} + \sum_{\alpha > 0}n(\alpha)H_{\alpha} -  2 \sum_{1 \leq i \leq a} \rho(H_{i})H_{i} + \sum_{1 \leq i \leq a} \rho(H_{i})^{2} - \sum_{\alpha > 0}n(\alpha) \rho(H_{\alpha}).$$ If we now observe, due to the orthonormality of the basis $\{H_{i} \}$ of $\mathfrak{a}$ relative to $B,$ that $ \sum_{1 \leq i \leq a} \rho(H_{i})H_{i}= H_{\rho} = \frac{1}{2} \sum_{\alpha > 0}n(\alpha)H_{\alpha}$ and $\sum_{1 \leq i \leq a} \rho(H_{i})^{2} = B(H_{\rho},H_{\rho})= \frac{1}{2} \sum_{\alpha > 0}n(\alpha) \rho(H_{\alpha}),$ we now have
$$\gamma(\omega)= \sum^{a}_{i=1} H^{2}_{i} + \sum_{\alpha > 0}n(\alpha)H_{\alpha} -  2 (\frac{1}{2} \sum_{\alpha > 0}n(\alpha)H_{\alpha}) + (\frac{1}{2} \sum_{\alpha > 0}n(\alpha) \rho(H_{\alpha})) - \sum_{\alpha > 0}n(\alpha) \rho(H_{\alpha})$$ $$\sum^{a}_{i=1} H^{2}_{i} - \frac{1}{2} \sum_{\alpha > 0}n(\alpha) \rho(H_{\alpha})= \sum^{a}_{i=1} H^{2}_{i}- B(H_{\rho},H_{\rho}).\;\Box $$

We can therefore get the expression for the radial component, $\delta^{'}(\omega),$ of $\omega$ using the above expression for $\gamma(\omega)$ in Theorem $2.9.$ This would give the most general form of the formula $$\delta^{'}(\omega^{'})=\frac{d^{2}}{dt^{2}}+ 2 \frac{\cosh 2t}{\sinh 2t} \frac{d}{dt}+1$$ of the case of $G=SL(2, \R)$ to all connected semisimple Lie group $G$ with finite center. The result is as follows.\\
\ \\

$\mathbf{3.4}$\textbf{ Proposition}($[2.],$ p. $133,$ and $[3(a.)],$ p. $269$). The radial component, $\delta^{'}(\omega),$ of the differential operator $\omega$ is given as $$\delta^{'}(\omega)= \sum_{1 \leq i \leq a}H^{2}_{i}+ 2 \sum_{\alpha > 0}n(\alpha)g_{\alpha}H_{\alpha} + \sum_{\alpha > 0}n(\alpha)H_{\alpha}.$$

\indent\textbf{\emph{Proof.}} Having known that $$\omega= - \sum^{m}_{k=1}Z^{2}_{k}+ \sum^{a}_{i=1} H^{2}_{i}+ \sum_{1 \leq i \leq n(\alpha), \alpha > 0}(X_{\alpha,i}X_{- \alpha,i}+ X_{- \alpha,i}X_{\alpha,i})$$ and that the first sum vanishes under $\beta_{\mathfrak{n}},$ it remains for us to find an expression for the sum $\sum_{1 \leq i \leq n(\alpha), \alpha > 0}(X_{\alpha,i}X_{- \alpha,i}+ X_{- \alpha,i}X_{\alpha,i})$ in terms of members, $f_{\alpha}$ and $g_{\alpha},$ of $\mathcal{R}^{+}_{0},$ as expected in the expression for $\delta^{'}(\omega)$ in Theorem $2.9$ and explicitly seen in the case of $G= SL(2,\R).$

To this end, let $X_{\alpha,i}= K_{\alpha,i}+ S_{\alpha,i}$ where $K_{\alpha,i} \in \mathfrak{t}$ and $S_{\alpha,i} \in \mathfrak{p}.$ Since $[X_{\alpha,i},X_{- \alpha,i}]= [K_{\alpha,i}+ S_{\alpha,i},K_{- \alpha,i}+ S_{- \alpha,i}]= [K_{\alpha,i}+ S_{\alpha,i},-K_{\alpha,i}+ S_{\alpha,i}]= 2 [K_{\alpha,i},S_{\alpha,i}] \in \mathfrak{p},$ and $[X,X^{'}] \equiv B(X,X^{'})H_{\alpha}(mod \mathfrak{m}), X \in \mathfrak{g}_{\alpha},X^{'} \in \mathfrak{g}_{- \alpha},$ we arrive at $[X_{\alpha,i},X_{- \alpha,i}]=H_{\alpha}.$ This then means that $[K_{\alpha,i},X{\alpha,i}]= [K_{\alpha,i},S_{\alpha,i}]= \frac{1}{2}[X_{\alpha,i},X_{- \alpha,i}]= \frac{1}{2}H_{\alpha}.$ Now setting $X=X_{\alpha,i}$ and $g=X_{\alpha,i}$ in $$\theta X \cdot g= f_{\alpha}(h)(X+ \theta X)^{h^{-1}}g-g_{\alpha}(h)(X+ \theta X)g$$ (from the equation $(*)$ of the proof of Theorem $2.7$) gives
\begin{align*}
2X_{- \alpha,i}X_{\alpha,i} &= -2 \theta(X_{\alpha,i})X_{\alpha,i}\\
&=  -2 [f_{\alpha}(h)(X_{\alpha,i}+ \theta X_{\alpha,i})^{h^{-1}}X_{\alpha,i}- g_{\alpha}(h)(X_{\alpha,i}+ \theta X_{\alpha,i})X_{\alpha,i}]  \\
&= -2 f_{\alpha}(h)[(K_{\alpha,i}+ S_{\alpha,i})+(K_{\alpha,i}- S_{\alpha,i})]^{h^{-1}}X_{\alpha,i}\\
&+ 2 g_{\alpha}(h)[(K_{\alpha,i}+ S_{\alpha,i})+(K_{\alpha,i}- S_{\alpha,i})]X_{\alpha,i}\\
&= -4 f_{\alpha}(h)K^{h^{-1}}_{\alpha,i}X_{\alpha,i}+ 4 g_{\alpha}(h)K_{\alpha,i}X_{\alpha,i}\\
&= -4 f_{\alpha}(h)K^{h^{-1}}_{\alpha,i}X_{\alpha,i}+ 4 g_{\alpha}(h)(X_{\alpha,i}K_{\alpha,i}+ [K_{\alpha,i},X_{\alpha,i}])\\
&= -4 f_{\alpha}(h)K^{h^{-1}}_{\alpha,i}X_{\alpha,i}+ 4 g_{\alpha}(h)X_{\alpha,i}K_{\alpha,i}+ 4g_{\alpha}(h)[K_{\alpha,i},X_{\alpha,i}]\\
&= -4 f_{\alpha}(h)K^{h^{-1}}_{\alpha,i}X_{\alpha,i}+ 4 g_{\alpha}(h)X_{\alpha,i}K_{\alpha,i}+ 4g_{\alpha}(h) ( \frac{1}{2}H_{\alpha})\\
&= -4 f_{\alpha}(h)K^{h^{-1}}_{\alpha,i}X_{\alpha,i}+ 4 g_{\alpha}(h)X_{\alpha,i}K_{\alpha,i}+ 2g_{\alpha}(h)H_{\alpha},
\end{align*}
so that
\begin{align*}
X_{\alpha,i}X_{- \alpha,i} + X_{- \alpha,i}X_{\alpha,i}\\
&= 2 X_{- \alpha,i}X_{\alpha,i} + H_{\alpha}\\
&= -4 f_{\alpha}(h)K^{h^{-1}}_{\alpha,i}X_{\alpha,i}+ 4 g_{\alpha}(h)X_{\alpha,i}K_{\alpha,i}+ 2 g_{\alpha}(h)H_{\alpha}+ H_{\alpha},
\end{align*}
which when substituted into the above expression for $\omega$ gives $$\omega= - \sum^{m}_{k=1}Z^{2}_{k}+ \sum^{a}_{i=1} H^{2}_{i} + \sum_{1 \leq i \leq n(\alpha)}(-4 f_{\alpha}(h)K^{h^{-1}}_{\alpha,i}X_{\alpha,i}+ 4 g_{\alpha}(h)X_{\alpha,i}K_{\alpha,i}+ 2 g_{\alpha}(h)H_{\alpha}+ H_{\alpha}).$$ Using the kernel of $\beta_{\mathfrak{n}},$ and hence of $\gamma_{\mathfrak{n}},$ in the expression for $\delta^{'}(\omega)$ we have
\begin{align*}
\delta^{'}(\omega)= -0+ \sum^{a}_{i=1} H^{2}_{i}-4(0)+4(0)+2 \sum_{\alpha > 0}n(\alpha)g_{\alpha}(h)H_{\alpha}+ \sum_{\alpha > 0}n(\alpha)H_{\alpha}
\end{align*}
$$= \sum^{a}_{i=1} H^{2}_{i}+ 2 \sum_{\alpha > 0}n(\alpha)g_{\alpha}(h)H_{\alpha}+ \sum_{\alpha > 0}n(\alpha)H_{\alpha}.\;\Box$$

As from now on we may start the discussion on the asymptotic behaviour of the spherical functions, $\varphi,$ from the pertubation theory of the system of differential equations satisfied by it given as $$\delta^{'}(\omega) \cdot \varphi_{\lambda} = \gamma_{\mathfrak{n}}(\lambda) \cdot \varphi_{\lambda},$$ with $\delta^{'}(\omega)$ as in Theorem $3.4.$ This is already contained in $[3(a.)],[3(b.)],[11.]$ and, more recently $[10.]$ In the present paper however, we shall seek to generalize the outlook that led, in $\S 2.,$ to the notion of a \textit{confluent spherical function} to all semisimple Lie groups, with real rank $1.$ This would require loading the structure of $G$ and its Lie algebra, $\mathfrak{g},$ into the known expressions for $\delta^{'}(\omega)$ and $\gamma_{\mathfrak{n}}(\lambda)$ as we do next.

We now take $G$ to be a connected semisimple Lie group with finite center and real rank $1.$ i.e., $dim(\mathfrak{a})=1.$ This implies that there exists exactly one simple root in $\Delta$ that we denote by $\alpha.$ This also means that $2 \alpha$ is the only other possible element in $\Delta^{+}$ and, if $p$ and $q$ represents the numbers $n(\alpha)$ and $n(2 \alpha)$ of roots in $\Delta$ which coincide on $\mathfrak{a}^{+}$ with $\alpha$ and $2 \alpha,$ respectively, then $$ \rho= \frac{1}{2} \sum_{\lambda \in \Delta^{+}} n(\lambda) \cdot \lambda= \frac{1}{2}(n(\alpha) \cdot \alpha + n(2 \alpha) \cdot (2 \alpha))= \frac{1}{2}(p+2q) \cdot \alpha$$ and that $p \geq 1, q \geq 0.$ Now choose $H_{0} \in \mathfrak{a}^{+}$ such that $\alpha(H_{0})=1.$ Since, for any $H \in \mathfrak{a},$ we always have $B(H,H)=2 \sum_{\lambda \in \Delta^{+}}n(\lambda) \cdot \lambda(H)^{2}=2(n(\alpha) \cdot (\alpha(H))^{2}+ n(2 \alpha) \cdot ((2 \alpha)(H))^{2})=2 \alpha(H)^{2}(p+4q),$ hence $$ B(H_{0},H_{0})=2 \cdot \alpha(H_{0})^{2}(p+4q)=2(p+4q) \cdots \cdots (i.),$$ and $B(H_{\alpha},H_{\alpha})=2 \cdot \alpha(H_{\alpha})^{2}(p+4q)).$ i.e., $\alpha(H_{\alpha})=2 \alpha(H_{\alpha})^{2}(p+4q).$ (since $\lambda(H^{'})=B(H_{\lambda},H^{'}),$ for all $H^{'} \in \mathfrak{a}.$) This implies that $\alpha(H_{\alpha})=(2(p+4q))^{-1}.$ Therefore $$H_{\alpha}= \alpha(H_{\alpha}) \cdot H_{0}=(2(p+4q))^{-1} \cdot H_{0} \cdots \cdots (ii.)$$ Equation $(ii.)$ and the relation $\rho= \frac{1}{2}(p+2q) \cdot \alpha$ also imply that $H_{\rho}= \rho(H_{\alpha}) \cdot H_{0}= \frac{1}{2}(p+2q) \cdot \alpha(H_{\alpha}) \cdot H_{0}= \frac{1}{2}(p+2q)(2(p+4q))^{-1} \cdot H_{0}.$ i.e., $$H_{\rho}= \frac{1}{2}(p+2q)(2(p+4q))^{-1} \cdot H_{0} \cdots \cdots (iii.)$$ We now state the major result of this section as follows. This result generalizes the situation of $SL(2,\R)$ above and motivates the concept of a confluent spherical function.\\

\noindent\textbf{{\Large \S 4. Confluent Spherical Functions And Their Algebra.}}

$\mathbf{4.1}$\textbf{ Theorem.} Let $G$ be a real rank $1$ connected semisimple Lie group with finite center and having the polar decomposition $G=K \cdot cl(A^{+}) \cdot K.$ Then every $K-$ biinvariant function on $G$ is spherical iff it is a hypergeometric function.

\indent\textbf{\emph{Proof.}} Using the isomorphism $t \longmapsto \exp tH_{0}, t \in \R$ between $\R$ and $A$ we identify $H_{0}$ with $\frac{d}{dt}$ so that
\begin{align*}
2(p+4q) \delta^{'}(\omega)\\
&=2(p+4q)[H^{2}_{\alpha}+2(n(\alpha)g_{\alpha}H_{\alpha}+n(2 \alpha)g_{2 \alpha}H_{2 \alpha})+(n(\alpha)H_{\alpha}+n(2 \alpha)H_{2 \alpha})]\\
&=2(p+4q)[(2(p+4q))^{-1} \cdot H^{2}_{0}+2(p \cdot \xi_{- \alpha}(\xi_{\alpha}-\xi_{- \alpha})^{-1} \cdot (2(p+4q))^{-1}H_{0}\\
&+2q \cdot \xi_{-2 \alpha}(\xi_{2 \alpha}- \xi_{-2 \alpha})^{-1}(2(p+4q))^{-1} \cdot H_{0})\\
&+(p \cdot (2(p+4q))^{-1} \cdot H_{0}+2q \cdot (2(p+4q))^{-1} \cdot H_{0})]\tag{from an adaptation of formula (ii.) above}\\
&=2(p+4q)(2(p+4q))^{-1}[\frac{d^{2}}{dt^{2}}+2(pg_{1} \frac{d}{dt}+2qg_{2} \frac{d}{dt})+(p \frac{d}{dt}+2q \frac{d}{dt})]\tag{where $g_{k}(t)=\xi_{-k \alpha}(t)(\xi_{k \alpha}(t)-\xi_{-k \alpha}(t))^{-1}=e^{-kt}(e^{kt}-e^{-kt})^{-1}=e^{-2kt}(1-e^{-2kt})^{-1},$ $k=1,2$ }\\
&= \frac{d^{2}}{dt^{2}}+ \{2(pg_{1}+2qg_{2})+(p+2q) \} \frac{d}{dt} \\
&= \frac{d^{2}}{dt^{2}}+ \{(2g_{1}+1)p+2(2g_{2}+1)q \} \frac{d}{dt}= \frac{d^{2}}{dt^{2}}+ \{p \coth t +2q \coth 2t \} \frac{d}{dt}\\
&= \frac{d^{2}}{dt^{2}}+ \{(p+q) \coth t +q \tanh t \} \frac{d}{dt} \cdots \cdots (iv) \tag{using the relation $\coth 2t = \frac{1}{2}(\coth t + \tanh t)$}\\
\end{align*}

In the case of $\gamma_{\mathfrak{n}}(\omega)(\lambda),$ we identify $\mathfrak{a}^{*}_{\mathbb{\C}}$ with $\mathbb{\C}$ via the map $\lambda \longmapsto \lambda(H_{0}),$ and for every $H_{1} \in \mathfrak{a},$ we set $B(H_{1},H_{1})=1.$ Now if we substitute $H=H_{1}$ into the relation $B(H,H)=2 \alpha(H)^{2}(p+4q),$ it gives $H_{1}=(2(p+4q))^{- \frac{1}{2}}H_{0},$ and when used with the expression for $H_{\rho}$ in $(iii.)$ above we have (from Lemma $3.3$) that
\begin{align*}
\gamma_{\mathfrak{n}}(\omega)(\lambda)= \lambda(H_{1})^{2}-B(H_{\rho},H_{\rho})\\
&= \lambda(H_{1})^{2}- \rho(H_{\rho})\\
&=(\lambda((2(p+4q))^{- \frac{1}{2}}H_{0}))^{2}- \rho(\frac{1}{2}(p+2q)(2(p+4q))^{-1}H_{0})\\
&= ((2(p+4q))^{-1} \lambda(H_{0}))^{2}- (\frac{1}{2}(p+2q)(2(p+4q))^{-1}) \rho(H_{0})\\
&=(2(p+4q))^{-1} \lambda^{2}-(\frac{1}{2}(p+2q)(2(p+4q))^{-1})(\frac{1}{2}(p+2q) \alpha(H_{0}))\\
&\tag{from the expression $\rho= \frac{1}{2}(p+2q) \alpha$}\\
\end{align*}
$$= (2(p+4q))^{-1}(\lambda^{2}- \frac{(p+2q)^{2}}{4}) \cdots \cdots (v.)$$

We now substitute the expressions for $\delta^{'}(\omega)$ and $\gamma_{\mathfrak{n}}(\omega)(\lambda)$ in $(iv.)$ and $(v.)$ into the equation $$\delta^{'}(\omega) \varphi_{\lambda}= \gamma_{\mathfrak{n}}(\omega)(\lambda) \cdot \varphi_{\lambda}$$ and define the function $f_{\lambda}$ on $\R,$ as $f_{\lambda}(t)= \varphi_{\lambda}(\exp tH_{0}),$ which is possible because of the above isomorphism between $\R$ and $A,$ to have $$(\frac{d^{2}}{dt^{2}}+ \{(p+q) \coth t +q \tanh t \} \frac{d}{dt}) f_{\lambda}= (\lambda^{2}- (\frac{(p+2q)^{2}}{4})) f_{\lambda} \cdots \cdots (\Upsilon \Upsilon).$$

This is the equation $(\Upsilon)$ at the tail-end of \S 2. for $G=SL(2, \R),$ where $p=2$ and $q=0.$ We now transform $(\Upsilon \Upsilon),$ as done in \S 2. for $\varphi^{''}_{\lambda}+2 \coth 2t \varphi^{'}_{\lambda}+(1- \lambda^{2}) \varphi_{\lambda}=0,$ by setting $z=-(\sinh t)^{2}.$ This implies that $\frac{dz}{dt}= -2 \sinh t \cosh t$ from which we may deduce that $\frac{d}{dt}=(-2 \sinh t \cosh t) \frac{d}{dz}$ and $\frac{d^{2}}{dt^{2}}=(4 \sinh^{2}t \cosh^{2}t) \frac{d^{2}}{dz^{2}}-2(\sinh^{2}t+ \cosh^{2}t) \frac{d}{dz}.$ Defining a function, $g_{\lambda},$ on $\mathbb{\C}$ as $g_{\lambda}(z)=f_{\lambda}(t)$ under the transformation $z=-(\sinh t)^{2}$ then converts equation $(\Upsilon \Upsilon)$ to $$(z(z-1) \frac{d^{2}}{dz^{2}}+((a+b+1)z-c) \frac{d}{dz}+ab)g_{\lambda}=0,$$ where the constants $a,b,c$ are given by $a= \frac{p+2q+2 \lambda}{4},$ $b= \frac{p+2q-2 \lambda}{4},$ $c= \frac{p+q+1}{2},$ respectively.

This is the well-known Gauss' hypergeometric equation. The point $z=0$ which corresponds to $t=0,$ is a regular singular point for this equation, and it is known that there is exactly one solution of it which is analytic at $z=0$ and takes the value $1$ there. This is the hypergeometric functions, $F(a,b,c:z),$ which, for $|z| < 1,$ is given as $$F(a,b,c:z)= \sum^{\infty}_{k=0} \frac{(a)_{k}(b)_{k}}{(c)_{k} \cdot k!}z^{k},$$ where $(m)_{k}:=m(m+1) \cdots (m+k-1).$ Now since $g_{\lambda}(0)=f_{\lambda}(0)=\varphi_{\lambda}(1)=1$ and $g_{\lambda}$ is analytic in $z$ at $z=0$ we conclude that $$\varphi_{\lambda}(\exp tH_{0})=F(a,b,c:-(\sinh t)^{2}),$$ where $t \in \R$ and $a,b,c$ are as given above.

Conversely, let a function $\varphi_{\lambda}$ be $K-$ biinvariant and be given as $\varphi_{\lambda}(\exp tH_{0})=F(a,b,c:-(\sinh t)^{2})$ for some $a,b,c.$ If we consider the equation $$\delta^{'}(\omega) \varphi_{\lambda}= \gamma_{\mathfrak{n}}(\omega)(\lambda) \cdot \varphi_{\lambda}$$ for some yet-to-be known constant $\gamma_{\mathfrak{n}}(\omega)(\lambda) \in \mathbb{\C},$ then $a= \frac{p+2q+2 \sqrt{\gamma_{\mathfrak{n}}(\omega)(\lambda)}}{4},$ $b= \frac{p+2q-2 \sqrt{\gamma_{\mathfrak{n}}(\omega)(\lambda)}}{4},$ and $c= \frac{p+q+1}{2}.$ For any known real rank $1$ connected semisimple Lie group $G,$ with finite center, in which $p$ and $q$ are also known, we may solve for $\gamma_{\mathfrak{n}}(\omega)(\lambda)$ explicitly from the above relations. With the fact that $\varphi_{\lambda}(1)=F(a,b,c:0)=1$ we conclude that $\varphi_{\lambda}$ is a spherical function on $G.\;\Box$

The above result shows the one-to-one correspondence between the hypergeometric functions and spherical functions on real rank $1$ semisimple Lie groups, $G.$ Now we recall the well-known notion of the \textit{confluent hypergeometric function} and use it, via Theorem $4.1$ above, to introduce the notion of a \textit{confluent spherical function} on $G,$ which is then later generalized using the \textit{Stanton-Tomas expansion} for spherical functions.

We recall that replacing $z(=-(\sinh t)^{2})$ by $\frac{z}{b}(= \frac{-(\sinh t)^{2}}{b})$ the hypergeometric equation gives $$z(z- \frac{z}{b}) \frac{d^{2}}{dz^{2}}g_{\lambda}(z)+ \{c-(1+ \frac{a+1}{b})z \} \frac{d}{dz}g_{\lambda}(z)-ag_{\lambda}(z)=0$$ becoming, as $b \longrightarrow \infty,$ $$z \frac{d^{2}}{dz^{2}}g_{\lambda}(z)+ (c-z) \frac{d}{dz}g_{\lambda}(z)-ag_{\lambda}(z)=0$$ whose solution, $g_{\lambda},$ is the confluent hypergeometric function, $_{1}F(a,c:z),$ is clearly given as $$_{1}F(a,c:z)= \lim_{b \longrightarrow \infty}F(a,b,c: \frac{z}{b})= \sum^{\infty}_{k=0} \frac{(a)_{k}}{(c)_{k}k!}z^{k}$$ where $a,c$ and $z$ are as above. Theorem $4.1$ implies that there exists a $K-$ biinvariant function, say $\varphi^{\sigma}_{\lambda},$ on $G$ such that $\varphi^{\sigma}_{\lambda}=\;_{1}F.$ It would be important to have a concise way of defining the function $\varphi^{\sigma}_{\lambda}.$ To do this we study more closely the properties of $_{1}F$ as follows.

The relationship between $z$ and $t,$ given as $z=-(\sinh t)^{2},$ could be recast as $t= sinh^{-1}(i \sqrt{z}).$ Now the process of deriving the confluent hypergeometric equation above entails substituting $z$ with $\frac{z}{b},$ before applying the limit as $b \longrightarrow \infty.$ Doing the same for the expression $t= sinh^{-1}(i \sqrt{z}),$ we have $t= sinh^{-1}(i \sqrt{\frac{z}{b}}).$ In applying the limit as $b \longrightarrow \infty,$ it follows that $t \longrightarrow 0.$ Now as $\lim_{b \longrightarrow \infty}F(a,b,c: \frac{z}{b})$ gives $_{1}F(a,c:z),$ the last statement above implies that we study $\lim_{t \longrightarrow 0} \varphi_{\lambda}(\exp tH_{0}).$ i.e., we study $\varphi_{\lambda}(\exp tH_{0})$ for sufficiently small values of $t.$ This observation is explicitly written as $$\varphi^{\sigma}_{\lambda}(\exp tH_{0})= \lim_{b \longrightarrow \infty} \varphi_{\lambda}(\exp tH_{0})_{\mid_{t= sinh^{-1}(i \sqrt{\frac{z}{b}})}}.$$

Since $b \longrightarrow \infty$ results to $t$ being very small, the equality above implies that the study of the confluent spherical functions, $\varphi^{\sigma}_{\lambda},$ on semisimple Lie groups is the same as the study of spherical functions, $\varphi_{\lambda}$ \textit{for sufficiently small values of t.} $i.e.,$ the study of the function $\varphi^{\sigma}_{\lambda}:G \rightarrow \C$ in which given $\epsilon > 0$ we can find $\delta=\delta(\epsilon) > 0$ such that $\mid \varphi_{\lambda}(\exp tH_{0}) - \varphi^{\sigma}_{\lambda}(\exp tH_{0})\mid < \epsilon$ whenever $t < \delta.$

This is the idea behind our notion of a confluent spherical function on $G,$ and to develop this idea further we make use of the \textit{Stanton-Tomas expansion} for spherical functions on a real rank $1$ connected semisimple Lie group, with finite center (See $[11.]$). We however start with a motivation via the case of $G=SL(2, \R)$ which proves the fact that \textit{Legendre functions admit a series expansion in terms of Bessel functions.}

It has been shown by Harish-Chandra, $[3(a.)],$ that every spherical function, on any connected semisimple Lie group $G,$ with finite center and arbitrary real rank, has the integral expansion $$\varphi_{\lambda}(x)= \int_{K}e^{(\lambda- \rho)(H(xk))}dk,$$ $\;\mbox{where}\;\lambda \in \mathfrak{a}^{*}_{\mathbb{\C}},\;\;x \in G,\;\;\rho= \frac{1}{2} \sum_{\alpha \in \Delta^{+}} dim(\mathfrak{g}_{\alpha}) \cdot \alpha.$ When $G=SL(2, \R),$ a calculation contained in $[15.],$ p. $339,$ shows that $$\varphi_{\lambda}(\exp tH_{0})= \frac{1}{2 \pi} \int^{2 \pi}_{0}(\cosh t+ \sinh t \cos \theta)^{(\lambda- \frac{1}{2})}d \theta.$$ This is the integral formula for the Legendre function, $P_{\lambda - \frac{1}{2}}(\cosh t)$ and is in consonance with the conclusion in \S 2. A change of contour in the integral yields $$\varphi_{\lambda}(\exp tH_{0})= c \int^{t}_{0} \cos (\lambda s)(\cosh t- \cosh s)^{- \frac{1}{2}}ds.$$ See $[11.]$ for some details. Now for small values of $t,$ it is known that $$(\cosh t- \cosh s)^{- \frac{1}{2}}=(t^{2}-s^{2})^{- \frac{1}{2}}+ error.$$ So that for sufficiently small values of $t,$ we have  $$\varphi_{\lambda}(\exp tH_{0})=c \int^{t}_{0} \cos (\lambda s)(t^{2}-s^{2})^{- \frac{1}{2}}ds= J_{0}(\lambda t),$$ where $J_{n}(\lambda t)$ is the Bessel function of order $n,$ giving as the series expansion $$J_{n}(x)= \frac{x^{n}}{2^{n} \cdot n!} \{1- \frac{x^{2}}{2^{2} \cdot 1! (n+1)}+ \frac{x^{4}}{2^{4} \cdot 2! (n+1)(n+2)}- \frac{x^{6}}{2^{6} \cdot 3!(n+1)(n+2)(n+3)}+ \cdots  \},$$ for any $n \in \mathbb{\Z},$ with $J_{-n}(x)=(-1)^{n}J_{n}(x),$ and, if $n \notin \mathbb{\Z},$ $$J_{n}(x)=\sum^{\infty}_{k=0}\frac{(-1)^{k}}{k! \Gamma(n+k+1)}(\frac{x}{2})^{n+2k}.$$  We can state our deductions above as follows.\\
\ \\

$\mathbf{4.2}$\textbf{ Theorem.} The confluent spherical functions, $\varphi^{\sigma}_{\lambda},$ on $G=SL(2, \R)$ are the zero-th order Bessel functions, $J_{0},$ of sufficiently small arguments.

\indent\textbf{\emph{Proof.}} Exactly as in the above deductions. $\Box$

This idea generalizes to all real rank $1$ semisimple Lie groups and is the first main result of Stanton and Tomas, $[11.].$ To state their result we make some preparations.

Let $n= dim(G/K),$ which is known to be equal to $p+q+1,$ define $c_{0}=c_{0}(G)$ and the function $\mathcal{J}_{\mu},$ respectively, as $$\pi^{\frac{1}{2}}2^{(q/2)-2} \frac{\Gamma(\frac{n-1}{2})}{\Gamma(\frac{n}{2})}$$ and $$\mathcal{J}_{\mu}(z) = \left\{\begin{array}{ll} \frac{J_{\mu}(z)}{z^{\mu}} \Gamma(\mu+ \frac{1}{2}) \Gamma(\frac{1}{2})2^{\mu -1}, & \mbox{if}\;
z \neq 0,\\
\;\;\;\;\;\;\;\;0, & \mbox{if}\; z =0,
\end{array}\right.$$
where $J_{\mu}$ is the Bessel function of order $\mu.$ Let also $D$ be the Jacobian for the polar decomposition of $G,$ then $D > 0$ on $A^{+}$ and is given by $$D(t)=D(a_{t})=e^{- \rho(\log a_{t})} \cdot \Pi_{\alpha > 0}(e^{2 \alpha(\log a_{t})}-1)^{n(\alpha)}.$$ (see \S 2. for its reduction in the case of $SL(2, \R)$)

In the case of a real rank $1$ group, $G,$ in which there are at most two positive roots, $\alpha$ and $2 \alpha,$ with multiplicity $p$ and $q,$ the Jacobian  reduces to $D(t)=e^{-2 \rho_{0}t}g_{1}(t)^{-p}g_{2}(t)^{-q}$ where $g_{k}(t)=e^{-2kt}(1-e^{-2kt})^{-1},\;k=1,2,$ and $\rho_{0}=\frac{1}{2}(p+2q).$ We now state a very important expansion formula for spherical functions, $\varphi_{\lambda},$ as follows.\\
\ \\

$\mathbf{4.3}$\textbf{ Theorem}($[11.],$ p. $253$) There exist $R_{0} >1, R_{1} >1,$ such that for any $t$ with, $0 \leq t \leq R_{0},$ the spherical function, $\varphi_{\lambda},$ has the given expansion $$\varphi_{\lambda}(\exp tH_{0})=c_{0}[\frac{t^{n-1}}{D(t)}]^{1/2} \sum^{\infty}_{m=0}t^{2m}a_{m}(t) \mathcal{J}_{\frac{(n-2)}{2}+m}(\lambda t)$$ where $$a_{0}(t)\equiv 1,\;\mbox{and} \mid a_{m}(t) \mid \leq cR^{-m}_{1}.\;\Box$$

The error on truncating the above series is controlled as in the following.\\
\ \\

$\mathbf{4.4}$\textbf{  Corollary}($[11.],$ p. $253$) There exist $R_{0} >1, R_{1} >1,$ such that for any $t$ with $0 \leq t \leq R_{0}$ and any $M \geq 0,$ the spherical function, $\varphi_{\lambda},$ is given as $$\varphi_{\lambda}(\exp tH_{0})=c_{0}[\frac{t^{n-1}}{D(t)}]^{1/2} \sum^{M}_{m=0}t^{2m}a_{m}(t) \mathcal{J}_{\frac{(n-2)}{2}+m}(\lambda t)+E_{M+1}(\lambda t)$$ where $$a_{0}(t)\equiv 1,\;\mbox{and} \mid a_{m}(t) \mid \leq cR^{-m}_{1},$$ and $$\mid E_{M+1}(\lambda t) \mid \leq c_{M}t^{2(M+1)},$$ if $\mid \lambda t \mid \leq 1,$ and $$\mid E_{M+1}(\lambda t) \mid \leq c_{M}t^{2(M+1)} \cdot (\lambda t)^{-(\frac{(n-1)}{2})+(M+1)},$$ if $\mid \lambda t \mid > 1.\;\Box$

We shall refer to the expansion in Theorem $4.3$ above as the \textit{Stanton-Tomas expansion} for spherical functions. We are therefore motivated to give the following general definition of a \textit{confluent spherical function} on a real rank $1$ semisimple Lie group $G.$\\
\ \\

$\mathbf{4.5}$\textbf{ Definition.} A \textit{confluent spherical function} is any $K-$ biinvariant function on $G(=K \cdot cl(A^{+}) \cdot K)$ which has the Stanton-Tomas expansion on $A^{+}.$

Explicitly, a function $\varphi^{\sigma}_{\lambda} \in C(G//K)$ is a confluent spherical function on $G$ if there exist $R_{0}> 1,\;\;R_{1}> 1$ such that for any $t,$ with $0\leq t\leq R_{0},$ and any $\lambda \in \C,$ $$\varphi^{\sigma}_{\lambda}(\exp tH_{0})=c_{0}[\frac{t^{n-1}}{D(t)}]^{1/2} \sum^{\infty}_{m=0}t^{2m}a_{m}(t) \mathcal{J}_{\frac{(n-2)}{2}+m}(\mid \lambda \mid t)$$ where $a_{0}(t)\equiv 1,\;\mbox{and} \mid a_{m}(t) \mid \leq cR^{-m}_{1}.$

We have $\varphi^{\sigma}_{0}(x)=0$ while the introduction of complex $\lambda$'s and their moduli guarantee that $\varphi^{\sigma}_{-\lambda}(x)=\varphi^{\sigma}_{\lambda}(x).$ The differential equation satisfied by the confluent spherical function is contained in the following.\\
\ \\

$\mathbf{4.6}$\textbf{ Theorem.} The function, $g_{\lambda}:\C \rightarrow \C$ which coincides with $\varphi^{\sigma}_{\lambda}$ on $K \times cl(A^{+}) \times K,$ via the transformation $z=-(\sinh t)^{2},$ for sufficiently small values of $t,$ is a solution of the differential equation $$z \frac{d^{2}y}{dz^{2}}+(c-z) \frac{dy}{dz}-ay=0.$$

\textbf{ Proof.} Since $z\frac{d^{2}}{dz^{2}}+(c-z)\frac{d}{dz}$ is a differential operator on $\C$ we may take it as the realization of some $q \in \mathfrak{Q},$ under the transformation $z=-(\sinh t)^{2},$ for sufficiently small values of $t.$ As $\varphi^{\sigma}_{\lambda}$ is, in particular, a spherical function on $G,$ with sufficiently small arguments, it satisfies the relation $q \cdot \varphi^{\sigma}_{\lambda}= \gamma(q)(\lambda) \varphi^{\sigma}_{\lambda}\;,$ for some $\gamma(q)(\lambda) \in \C.$  The result follows if we set $\gamma(q)(\lambda)=a.\;\Box$\\
\ \\

$\mathbf{4.7}$\textbf{ Remarks.}

The confluent spherical functions, $\varphi^{\sigma}_{\lambda}(\exp t H_{0}),$ are defined for all $t \geq 0,\;\mbox{and all}\;\lambda \in \C.$ Indeed, $\varphi^{\sigma}_{\lambda}(1)=\varphi^{\sigma}_{\lambda}(\exp (0) H_{0})=0,$ since $D > 0$ on $A^{+}.$

We denote the set of all confluent spherical functions on $G$ by $C^{\sigma}(G)$ and consider it as an algebra in the following precise manner:\\
\ \\

$\mathbf{4.8}$\textbf{ Definition.} A non-empty set $\mathbb{A},$ whose entries are indexed by a set $\Delta,$ is called a \textit{$\Delta-$algebra} (over a field $\mathbb{K}$) if it is an (associative) algebra with respect to the operations
\begin{description}
\item [{(i.)}] $a_{\lambda_{1}}+a_{\lambda_{2}}:=a_{\lambda_{1}+\lambda_{2}},$ for every $a_{\lambda_{1}},\;a_{\lambda_{2}} \in \mathbb{A}.$
\item [{(ii.)}] $\alpha a_{\lambda}:=a_{\alpha \lambda},$ for every $a_{\lambda} \in \mathbb{A},\;\alpha \in \mathbb{K}.$
\item [{(iii.)}] $a_{\lambda_{1}} \cdot a_{\lambda_{2}}:=a_{\lambda_{1} \lambda_{2}},$ for every $a_{\lambda_{1}},\;a_{\lambda_{2}} \in \mathbb{A}.$
\end{description}
We shall refer to $(i.)-(iii.)$ above as the \textit{$\Delta-$operations} on $\mathbb{A}.$

One of our major results in this seminar is the following Theorem. In order to establish this result we denote the Schwartz algebra of spherical functions on $G$ by $\mathcal{C}(G//K)$ and equip both $C^{\sigma}(G)$ and $\mathcal{C}(G//K)$ with the \textit{$\mathfrak{a}^{*}_{\C}-$operations.}

Let $\mathfrak{w}$ denote the Weyl group of $(\mathfrak{g},\mathfrak{a})$ and $$l_{s}(\mathfrak{a}^{*}_{\C})=\{\lambda \in \mathfrak{a}^{*}_{\C}:s^{-1}\lambda=\lambda\} $$ for some $s \in \mathfrak{w}.$ Clearly $l_{s^{-1}}(\mathfrak{a}^{*}_{\C})=l_{s}(\mathfrak{a}^{*}_{\C}),\;\forall\;s \in \mathfrak{w}$ and, if $id$ represent the identity element of $\mathfrak{w},$ then $l_{id}(\mathfrak{a}^{*}_{\C})=\mathfrak{a}^{*}_{\C}.$ In general, $l_{s}(\mathfrak{a}^{*}_{\C})\subseteq\mathfrak{a}^{*}_{\C},\;\forall\;s \in \mathfrak{w}.$

We shall refer to a map between any two $\Delta-$algebras as being \textit{$\Delta-$linear} if it preserves $(i.)$ and $(ii.)$ of Definition $4.8.$\\
\ \\

$\mathbf{4.9}$\textbf{ Theorem.} The sets $\mathcal{C}(G//K)$ and $C^{\sigma}(G)$ are $\mathfrak{a}^{*}_{\C}-$algebras over $\C,$ where the zero and identity elements of $\mathcal{C}(G//K)$ are $\Xi$ and $\varphi_{1},$ respectively. The map $\sigma:\mathcal{C}(G//K)\longrightarrow C^{\sigma}(G),$ given by $\sigma(\varphi_{\lambda})=\varphi^{\sigma}_{\lambda}$ is non-trivial and well-defined up to $l_{s}(\mathfrak{a}^{*}_{\C}),\;s \in \mathfrak{w}.$ Moreover $\sigma$ is an $\mathfrak{a}^{*}_{\C}-$linear map and an isomorphism for all real-positive $\lambda.$

\indent\textbf{\emph{Proof.}} We verify using the $\mathfrak{a}^{*}_{\C}-$operations that $\mathcal{C}(G//K)$ is a $\mathfrak{a}^{*}_{\C}-$algebra. The situation for $C^{\sigma}(G)$ follows the same pattern.\\
To this end let $\varphi_{\lambda_{1}},\;\varphi_{\lambda_{2}},\;\varphi_{\lambda_{3}},\;\varphi_{\lambda} \in \mathcal{C}(G//K),\;\mbox{and}\;\alpha,\;\beta\in \C;$\\
$(i.)$ $\varphi_{\lambda_{1}}+\varphi_{\lambda_{2}}=\varphi_{\lambda_{1}+\lambda_{2}} \in \mathcal{C}(G//K),$ since $\lambda_{1}+\lambda_{2} \in \mathfrak{a}^{*}_{\C}.$\\
$(ii.)$ It is also clear that $\varphi_{\lambda_{1}}+(\varphi_{\lambda_{2}}+\varphi_{\lambda_{3}})=(\varphi_{\lambda_{1}}+\varphi_{\lambda_{2}})+\varphi_{\lambda_{3}}.$\\
$(iii.)\;\Xi+\varphi_{\lambda}=\varphi_{0}+\varphi_{\lambda}=\varphi_{0+\lambda}=\varphi_{\lambda}=\varphi_{\lambda+0}=\varphi_{\lambda}+\varphi_{0}
=\varphi_{\lambda}+\Xi.$\\
$(iv.)$ $\varphi_{-\lambda}+\varphi_{\lambda}=\varphi_{-\lambda+\lambda}=\varphi_{0}=\varphi_{\lambda+(-\lambda)}=\varphi_{\lambda}+\varphi_{-\lambda}.$\\
$(v.)\;\varphi_{\lambda_{1}}+\varphi_{\lambda_{2}}=\varphi_{\lambda_{1}+\lambda_{2}}=\varphi_{\lambda_{2}+\lambda_{1}}=
\varphi_{\lambda_{2}}+\varphi_{\lambda_{1}}$
$(vi.)\;\alpha\varphi_{\lambda}=\varphi_{\alpha\lambda}\in \mathcal{C}(G//K),$ since $\alpha\lambda \in \mathfrak{a}^{*}_{\C}.$\\
$(vii.)\;\alpha(\beta\varphi_{\lambda})=\varphi_{\alpha(\beta\lambda)}=\varphi_{(\alpha\beta)\lambda}=(\alpha\beta)\varphi_{\lambda}.$\\
$(viii.)\;1\varphi_{\lambda}=\varphi_{1\lambda}=\varphi_{\lambda}=\varphi_{\lambda1}=\varphi_{\lambda}1.$\\
$(ix.)\;\alpha(\varphi_{\lambda_{1}}+\varphi_{\lambda_{2}})=\varphi_{\alpha(\lambda_{1}+\lambda_{2})}=\varphi_{\alpha\lambda_{1}+\alpha\lambda_{2}}
=\alpha\varphi_{\lambda_{1}}+\alpha\varphi_{\lambda_{2}}.$\\
$(x.)\;(\alpha+\beta)\varphi_{\lambda}=\varphi_{(\alpha+\beta)\lambda}=\varphi_{\alpha\lambda+\beta\lambda}
=\alpha\varphi_{\lambda}+\beta\varphi_{\lambda}.$\\
$(xi.)\;\varphi_{\lambda_{1}}\cdot\varphi_{\lambda_{2}}=\varphi_{\lambda_{1}\cdot\lambda_{2}} \in \mathcal{C}(G//K),$ since $\lambda_{1}\cdot\lambda_{2} \in \mathfrak{a}^{*}_{\C}.$\\
$(xii.)$ It is also clear that $\varphi_{\lambda_{1}}\cdot(\varphi_{\lambda_{2}}\cdot\varphi_{\lambda_{3}})
=(\varphi_{\lambda_{1}}\cdot\varphi_{\lambda_{2}})\cdot\varphi_{\lambda_{3}}.$\\
$(xiii.)\;\varphi_{1}\cdot\varphi_{\lambda}=\varphi_{1\cdot\lambda}=\varphi_{\lambda}=\varphi_{\lambda\cdot1}=\varphi_{\lambda}\cdot\varphi_{1}.$\\
$(xiv)\;(\varphi_{\lambda_{1}}+\varphi_{\lambda_{2}})\cdot\varphi_{\lambda_{3}}=\varphi_{(\lambda_{1}+\lambda_{2})\cdot\lambda_{3}}
=\varphi_{(\lambda_{1}\cdot\lambda_{3})+(\lambda_{2}\cdot\lambda_{3})}
=(\varphi_{\lambda_{1}}\cdot\varphi_{\lambda_{3}})+(\varphi_{\lambda_{2}}\cdot\varphi_{\lambda_{3}}).$\\
$(xv.)\;\varphi_{\lambda_{1}}\cdot(\varphi_{\lambda_{2}}+\varphi_{\lambda_{3}})=\varphi_{\lambda_{1}\cdot(\lambda_{2}+\lambda_{3})}
=\varphi_{(\lambda_{1}\cdot\lambda_{2})+(\lambda_{1}\cdot\lambda_{3})}
=(\varphi_{\lambda_{1}}\cdot\varphi_{\lambda_{2}})+(\varphi_{\lambda_{1}}\cdot\varphi_{\lambda_{3}}),$ verifying the first statement.

The map $\sigma$ is \underline{non-trivial} from Theorem $4.3$ and Definition $4.5.$

We observe also that $\varphi_{\lambda_{1}}=\varphi_{\lambda_{2}}$ iff $\lambda_{2}=s\lambda_{1},$ for some $s \in \mathfrak{w};\;[2.],\;p.\;106.$ Since $s^{-1}\lambda_{2}=\lambda_{2}$ iff $s=id \in \mathfrak{w},$ it follows therefore $\sigma(\varphi_{\lambda_{1}})=\sigma(\varphi_{s^{-1}\lambda_{2}})=\sigma(\varphi_{\lambda_{2}}).$ Hence $\varphi_{\lambda_{1}}=\varphi_{\lambda_{2}}$ implies $\sigma(\varphi_{\lambda_{1}})=\sigma(\varphi_{\lambda_{2}}),$ showing that $\sigma$ is \underline{well-defined} up to $l_{id}(\mathfrak{a}^{*}_{\C})=\mathfrak{a}^{*}_{\C}.$\\

$\sigma$ is a \underline{$\mathfrak{a}^{*}_{\C}-$linear map} since $\sigma(\varphi_{\lambda_{1}}+\varphi_{\lambda_{2}})=\sigma(\varphi_{\lambda_{1}+\lambda_{2}})
=\varphi^{\sigma}_{\lambda_{1}+\lambda_{2}}=\varphi^{\sigma}_{\lambda_{1}}+\varphi^{\sigma}_{\lambda_{2}}
=\sigma(\varphi_{\lambda_{1}})+\sigma(\varphi_{\lambda_{2}})$\\ and,
$\sigma(\alpha\varphi_{\lambda})=\sigma(\varphi_{\alpha\lambda})=\varphi^{\sigma}_{\alpha\lambda}
=\alpha\varphi^{\sigma}_{\lambda}=\alpha\sigma(\varphi_{\lambda})$.\\
\underline{onto}:\\
Let $\varphi^{\sigma}_{\lambda}\in C^{\sigma}(G),$ then $\lambda \in l_{s}(\mathfrak{a}^{*}_{\C}),\;s \in \mathfrak{w}.$ Therefore $\varphi_{\lambda} \in \mathcal{C}(G//K).$ Hence $\sigma(\varphi_{\lambda})=\varphi^{\sigma}_{\lambda}.$\\
\underline{into}:\\
Let $\sigma(\varphi_{\lambda_{1}})=\sigma(\varphi_{\lambda_{2}}),$ then from Theorem $4.3,$ $(\lambda_{1} t)^{n+2k-m}=(\lambda_{2} t)^{n+2k-m}.$ Hence $$(n+2k-m)\log\lambda_{1}=(n+2k-m)\log\lambda_{2},\;t\neq0.$$ We then have $\lambda_{1}=\lambda_{2},$ which makes sense only if $\lambda_{1},\;\lambda_{2}\in \R^{+}.\;\;\Box$\\

We shall consider the extension of these results to a connected semisimple Lie group with finite center and of real rank $m>0$ in another paper.
\ \\

{\bf   References.}
\begin{description}
\item [{[1.]}] Barchini, L. and Zierau, R., \textit{Differential operators on homogeneous spaces,} Notes of Lectures at the ICE-EM Australian Graduate School in Mathematics in Brisbane, $2$ nd- $20$ th July, $(2007).$ \textit{Retrieved from http:// www.math.okstate.edu/zierau/papers.
    html} on 12 th January, $2011.$
       \item [{[2.]}] Gangolli, R. and Varadarajan, V. S., \textit{Harmonic analysis of spherical functions on real reductive groups,} Ergenbnisse der Mathematik und ihrer Grenzgebiete, \textbf{101}, Springer-Verlag, $1988.$
            \item [{[3.]}] Harish-Chandra, (a.) Spherical functions on a semisimple Lie group I, \textit{Amer. J. Math.,} \textbf{80}, $(1958),$ pp. $241-310.$ (b.) Spherical functions on a semisimple Lie group II, \textit{Amer. J. Math.,} \textbf{80}, $(1958),$ pp. $553-613.$
            \item [{[4.]}] Helgason, S., \textit{Differential geometry and symmetric spaces,} Academic Press, $1962.$
                \item [{[5.]}] Hirano, M., Ishii, T., and Oda, T., confluent from Siegel-Whittaker functions to Whittaker functions on $Sp(2, \R),$ ($2004$), pp. $1$-$14.$
                    \item [{[6.]}] Knapp, A. W., \textit{Lie groups beyond an introduction,} Progress in Mathematics, \textbf{140}, Birk$\ddot{a}$user Verlag, $2002.$
                    \item [{[7.]}] Knapp, A. W., and Trapa, P. E., \textit{Representations of semisimple Lie groups,} IAS/Park City Mathematics Series, \textbf{8}, ($1998$), pp. $1$-$83.$
                    \item [{[8.]}] Lang, S., $SL_{2}(\R),$ Addison-Wesley Publishing Company, $1975.$
                \item [{[9.]}] Oyadare, O. O., Harmonic analysis of spherical functions on $SL(2, \R),$ \textit{Unpublished M.Sc. Dissertation,} University of Ibadan, $2007.$
                \item [{[10.]}] Sawyer, P., The asymptotic expansion of spherical functions on symmetric cones, \textit{Pacific J. Math.,} \textbf{200}($1$), $(2001),$ pp. $251$-$256.$
                    \item [{[11.]}] Stanton, R. J., and Tomas, P. A., Expansion for spherical functions on noncompact symmetric spaces, \textit{Acta. Math.,} \textbf{140}, ($1978$), pp. $251$-$276.$
                        \item [{[12.]}] Varadarajan, V. S., \textit{Lie groups, Lie algebras and their representations,} Springer-Verlag, $1984.$
                \item [{[13.]}] Varadarajan, V. S., \textit{An introduction to harmonic analysis on semisimple Lie groups,} Cambridge University Press, Cambridge, $1989.$
                    \item [{[14.]}] Warner, F., \textit{Foundations of differentiable manifolds and Lie groups,} Scott, Foreman and Company, Glenview, $1971.$
                \item [{[15.]}] Warner, G., \textit{Harmonic analysis on semisimple Lie groups, II}, Grundlehren der mathematischen Wissenschaften, Band \textbf{189},Springer-Verlag, $1972.$
                    \item [{[16.]}] Whittaker, E. T., and Watson, G. N., \textit{A course of modern analysis,} 4 th edition, Cambridge University Press, Cambridge, $1965.$
                        \end{description}

\end{document}